\documentclass[11pt,a4paper]{article}

\usepackage{soul}
\usepackage{tikz}
\usepackage{amsthm}
\usepackage{amsmath}
\usepackage{amssymb}
\usepackage{mathrsfs}

\usepackage{geometry}
\usepackage{graphicx}
\usepackage{amsfonts}
\usepackage{epstopdf}
\usepackage{placeins}
\usepackage{enumerate}
\usepackage{enumitem}
\usepackage{subcaption}
\usepackage{microtype}
\usepackage[hidelinks]{hyperref} 
\usepackage[normalem]{ulem} 

\usetikzlibrary{calc,arrows.meta, positioning,decorations.markings}

\numberwithin{equation}{section}

\allowdisplaybreaks

\newcommand{\Rmnum}[1]{\uppercase\expandafter{\romannumeral#1}} 

\def\Xint#1{\mathchoice
{\XXint\displaystyle\textstyle{#1}}%
{\XXint\textstyle\scriptstyle{#1}}%
{\XXint\scriptstyle\scriptscriptstyle{#1}}%
{\XXint\scriptscriptstyle\scriptscriptstyle{#1}}%
\!\int}
\def\XXint#1#2#3{{\setbox0=\hbox{$#1{#2#3}{\int}$ }
\vcenter{\hbox{$#2#3$ }}\kern-.6\wd0}}

\def\dashint{\Xint-}

\theoremstyle{plain}
\newtheorem{theorem}{Theorem}[section]
\newtheorem{proposition}[theorem]{Proposition}
\newtheorem{lemma}[theorem]{Lemma}
\newtheorem{corollary}[theorem]{Corollary}
\newtheorem{definition}[theorem]{Definition}

\theoremstyle{definition}
\newtheorem{remark}[theorem]{Remark}

\renewcommand{\thefootnote}{}

\AtEndDocument{
\bigskip{\footnotesize%
   \textsc{Department of Mathematics, Aarhus University, 8000 Aarhus C, Denmark} \par
  \textit{E-mail address}: \texttt{yangmengqh@gmail.com}\par
  \textit{E-mail address}: \texttt{yang@math.au.dk}
}}

\begin{document}

\title{\Large{On the dichotomy of $p$-walk dimensions on metric measure spaces II}}
\author{Meng Yang}
\date{}

\maketitle

\abstract{On a volume doubling metric measure space endowed with a family of $p$-energies such that the Poincar\'e inequality and the capacity upper bound with $p$-walk dimension $\beta_p$ hold, for $p$ in an open interval $I\subseteq(1,+\infty)$, we prove the following dichotomy: either $\beta_p=p$ for \emph{all} $p\in I$, or the function $p\mapsto \frac{\beta_p}{p}$ takes values in $(1,+\infty)$ and is \emph{strictly} decreasing on $I$. Moreover, if the metric measure space is $d_h$-Ahlfors regular, we prove the following dichotomy: either $\beta_p=d_h+p-1$ for \emph{all} $p\in I$, or the function $p\mapsto \frac{\beta_p-d_h}{p-1}$ takes values in $(-\infty,1)$ and is \emph{strictly} increasing on $I$.}

\footnote{\textsl{Date}: \today}
\footnote{The author is grateful to Joseph Feneuil for sharing, in February 2025, a preprint version of \cite{Fen26} on Riesz transforms on fractal-like spaces, which motivated \cite{Yan25e} and the study of the monotonicity of $p\mapsto \frac{\beta_p}{p}$; see Remark \ref{rmk_Fen} for a discussion of this connection. The author is also grateful to Naotaka Kajino and Ryosuke Shimizu for helpful discussions on the monotonicity of a related exponent arising in their work \cite{KS24sing} on nested fractals. This exponent corresponds to $\frac{\beta_p-d_h}{p-1}$ in the notation of this paper.}

\renewcommand{\thefootnote}{\arabic{footnote}}
\setcounter{footnote}{0}

\section{Introduction}\label{sec_intro}

Diffusions on fractals often exhibit scaling behavior different from that of Brownian motion in Euclidean spaces. This behavior is described by the walk dimension $\beta_2$ appearing in sub-Gaussian heat kernel estimates. For example, on the standard Sierpi\'nski gasket, $\beta_2=\log 5/\log 2>2$ (see \cite{BP88,Kig93}), whereas on the standard Sierpi\'nski carpet, $\beta_2\approx2.09697$ (see \cite{BB89,BBS90,KZ92,BB92}). The corresponding local regular Dirichlet forms provide analogues of the classical Dirichlet integral $\int_{\mathbb{R}^d}|\nabla f(x)|^2\mathrm{d}x$; see \cite{FOT11}.

For $p\in(1,+\infty)\backslash\{2\}$, the non-linear counterpart is a $p$-energy, which generalizes
$$
\int_{\mathbb{R}^d}|\nabla f(x)|^p\mathrm{d}x.
$$
Following the early work in \cite{HPS04}, the study of $p$-energies on fractals and general metric measure spaces has advanced considerably. Recent developments include constructions and analysis on p.c.f. self-similar sets \cite{CGQ22}, Sierpi\'nski carpets \cite{Shi24,MS25}, Vicsek-type fractals \cite{BC23,CGYZ26}, conductively homogeneous compact metric spaces \cite{Kig23}, Cheeger spaces \cite{AB25}, and Laakso-type fractal spaces \cite{AEBS25}. Further developments in non-linear potential theory and a survey of the subject can be found in \cite{KS24aV5,KS27}. In this setting, the $p$-walk dimension $\beta_p$ describes the spatial scaling in the associated functional inequalities, including the Poincar\'e inequality, and in the associated capacity estimates. For $p=2$, $\beta_2$ agrees with the walk dimension of the corresponding diffusion.

Continuing the study initiated in \cite{Yan25e}, we investigate the dependence of $\beta_p$ on $p$. On the Vicsek set, the $p$-walk dimension is explicitly given by $\beta_p=p+d_h-1$, where $d_h=\log 5/\log 3$ is the Hausdorff dimension; see \cite{BC23}. On the Sierpi\'nski gasket and the Sierpi\'nski carpet, $\beta_p>p$ for any $p\in(1,+\infty)$; see \cite[Theorems 9.8 and 9.13]{KS24aV5}. These examples motivate the study of the regularity and monotonicity of $\beta_p$ as a function of $p$ on general metric measure spaces.

The Besov characterization of $p$-energies in \cite{Shi24a}, together with the results on critical Besov exponents in \cite{Bau24}, shows, under volume doubling and suitable functional inequalities, that $p\mapsto\beta_p$ is continuous and non-decreasing, while $p\mapsto \frac{\beta_p}{p}$ is non-increasing. In \cite{Yan25e}, we established a dichotomy: under volume doubling and suitable functional inequalities for a family of $p$-energies indexed by $p$ in an interval $I\subseteq(1,+\infty)$, either $\beta_p=p$ for any $p\in I$, or $\beta_p>p$ for any $p\in I$. The present paper supersedes that preprint. We strengthen the second alternative by proving that $p\mapsto \frac{\beta_p}{p}$ is strictly decreasing on $I$. This strict monotonicity is indeed related to the unboundedness of Riesz transforms on fractal-like spaces studied in \cite{Fen26}; see Remark \ref{rmk_Fen}.

On a $d_h$-Ahlfors regular space satisfying the chain condition, the Poincar\'e inequality and the capacity upper bound with exponent $\beta_p$ imply
$$
p\le\beta_p\le d_h+p-1;
$$
see \cite[Proposition~2.1]{Yan25a} and Figure \ref{fig_betapdh}. This leads to the normalized exponent
$$
\gamma_p:=\frac{\beta_p-d_h}{p-1}\le1.
$$
The regularity and monotonicity of $p\mapsto\gamma_p$ can provide information about the analytic structure of the space. Roughly speaking, when $\beta_p>d_h$, this quantity arises as a H\"older exponent for $p$-harmonic functions in the resistance form setting; see \cite[Theorem 6.31]{KS24aV5}. The same exponent also describes the relation between the resistance metric and the original metric. Indeed, under the above conditions and when $\beta_p>d_h$, the function
$$
\rho^{(p)}(x,y):=R^{(p)}(x,y)^{\frac{1}{p-1}}\asymp d(x,y)^{\gamma_p}
$$
is a metric, called the $p$-resistance metric; see \cite[Corollary 6.32 and Definition 6.33]{KS24aV5}.

Strict monotonicity of $\gamma_p$ is closely connected with the singularity of energy measures associated with distinct values of $p$. On the standard two-dimensional Sierpi\'nski gasket, it is proved in \cite[Theorems 4.8 and 4.9]{KS27} that $p\mapsto\gamma_p$ is strictly increasing, and the resulting strict inequality implies mutual singularity of $p$- and $q$-energy measures for $p\ne q$. In \cite{KS24sing}, they obtain corresponding strict monotonicity and singularity results for a class of affine nested fractals, including higher-dimensional and higher-level Sierpi\'nski gaskets; see also the announcement in \cite[p. 14]{KS24aV5}.

For the golden ratio Sierpi\'nski gasket, non-linear renormalization was studied and $p$-resistance forms were constructed in \cite{ChenQiu26}. It was proved that the standard renormalization factor $r_p\in(0,1)$ is continuous and strictly decreasing in $p\in(1,+\infty)$, while $p\mapsto r_p^{\frac{1}{p-1}}$ is continuous and non-increasing; see \cite[Theorem 3]{ChenQiu26}. Their strict monotonicity result concerns the renormalization factor itself.

These results, together with the identity $\gamma_p=1$ on the Vicsek set, suggest a second dichotomy on general metric measure spaces: either $\gamma_p=1$ throughout the interval $I$, or $\gamma_p<1$ throughout the interval $I$ and $p\mapsto\gamma_p$ is strictly increasing on $I$. We prove this dichotomy. Both dichotomies in this paper follow from geometric and analytic assumptions on the metric measure space and its $p$-energies, without requiring a self-similar structure.

\section{Statement of main results}\label{sec_result}

Before stating our main results, we fix some notation. Throughout this paper, $C$, $C_1$, $C_2$, $C_A$, $C_B$ denote positive constants whose values may change at each occurrence. The sign $\asymp$ means that the ratio of the two sides is bounded from above and below by positive constants. The sign $\lesssim$ ($\gtrsim$) means that the LHS is bounded above (below) by a positive constant times the RHS. We use $x_+$ to denote the positive part of $x\in \mathbb{R}$, that is, $x_+=\max\{x,0\}$. For two Borel measures $\mu, \nu$, we write $\mu\le\nu$ if $\mu(A)\le \nu(A)$ for any Borel set $A$, and write $\mu\lesssim\nu$ if $\mu\le C\nu$ for some constant $C>0$. When $\mu,\nu$ are $\sigma$-finite, $\mu\le \nu$ is equivalent to $\mu\ll \nu$ and $\frac{\mathrm{d}\mu}{\mathrm{d}\nu}\le1$ $\nu$-a.e. We use $\#A$ to denote the cardinality of a set $A$.

Let $(X,d,m)$ be a complete unbounded metric measure space, that is, $(X,d)$ is a complete unbounded locally compact separable metric space and $m$ is a positive Radon measure on $X$ with full support. Throughout this paper, we always assume that all metric balls are relatively compact. For any $x\in X$ and any $r>0$, denote $B(x,r)=\{y\in X:d(x,y)<r\}$ and $V(x,r)=m(B(x,r))$. If $B=B(x,r)$, then denote $\delta B=B(x,\delta r)$ for any $\delta>0$. Let $\mathcal{B}(X)$ be the family of all Borel measurable subsets of $X$. Let $C(X)$ be the family of all continuous functions on $X$. Let $C_c(X)$ be the family of all continuous functions on $X$ with compact support. Denote $\dashint_A=\frac{1}{m(A)}\int_A$ and $u_A=\dashint_Au\mathrm{d} m$ for any measurable set $A$ with $m(A)\in(0,+\infty)$ and any function $u$ such that the integral $\int_Au\mathrm{d} m$ is well-defined.


Let $\varepsilon>0$. We say that $W\subseteq X$ is an $\varepsilon$-net if for any distinct $x,y\in W$, we have $d(x,y)\ge \varepsilon$, and for any $z\in X$, there exists $x\in W$ such that $d(x,z)<\varepsilon$. Since $(X,d)$ is separable, all $\varepsilon$-nets are countable.

We say that the chain condition \ref{eq_CC} holds if there exists a constant $C_{cc}>0$ such that for any $x,y\in X$ and any integer $n\ge1$, there exists a sequence $\{x_k:0\le k\le n\}$ of points in $X$ with $x_0=x$ and $x_n=y$ such that
\begin{equation*}\label{eq_CC}\tag*{CC}
d(x_k,x_{k-1})\le C_{cc} \frac{d(x,y)}{n}\text{ for any }k=1,\ldots,n.
\end{equation*}
Throughout this paper, we always assume \ref{eq_CC}.

We say that the volume doubling condition \ref{eq_VD} holds if there exists a constant $C_{VD}>0$ such that
\begin{equation*}\label{eq_VD}\tag*{VD}
V(x,2r)\le C_{VD}V(x,r)\text{ for any }x\in X,\ r>0.
\end{equation*}
Under \ref{eq_VD}, the metric space $(X,d)$ is doubling, that is, there exists an integer $N\ge1$, depending only on $C_{VD}$, such that for any $r>0$, any ball in $(X,d)$ with radius $2r$ can be covered by at most $N$ balls with radius $r$; see \cite[Remarque on p. 67]{CW71} or \cite[pp. 81--82]{Hei01}.

For $d_h>0$, we say that the Ahlfors regular condition \ref{eq_Vdh} holds if there exists a constant $C_{R}>0$ such that
\begin{equation*}\label{eq_Vdh}\tag*{$\text{V}(d_h)$}
\frac{1}{C_{R}}r^{d_h}\le V(x,r)\le C_{R}r^{d_h}\text{ for any }x\in X,\ r>0.
\end{equation*}

We say that $(\mathcal{E},\mathcal{F})$ is a $p$-energy on $(X,d,m)$ if $\mathcal{F}$ is a dense subspace of $L^p(X;m)$ and $\mathcal{E}:\mathcal{F}\to[0,+\infty) $ satisfies the following conditions.

\begin{enumerate}[label=(\arabic*)]
\item $\mathcal{E}^{1/p}$ is a seminorm on $\mathcal{F}$, that is, for any $f,g\in\mathcal{F}$, $c\in\mathbb{R}$, we have $\mathcal{E}(f)\ge0$, $\mathcal{E}(cf)^{1/p}=|c|\mathcal{E}(f)^{1/p}$, and $\mathcal{E}(f+g)^{1/p}\le\mathcal{E}(f)^{1/p}+\mathcal{E}(g)^{1/p}$.
\item (Closed property) $(\mathcal{F},(\mathcal{E}(\cdot)+\lVert {\cdot}\rVert^p_{L^p(X;m)})^{1/p})$ is a Banach space.
\item (Markovian property) For any $\varphi\in C(\mathbb{R})$ with $\varphi(0)=0$ and $|\varphi(t)-\varphi(s)|\le|t-s|$ for any $t, s\in\mathbb{R}$, for any $f\in\mathcal{F}$, we have $\varphi(f)\in\mathcal{F}$ and $\mathcal{E}(\varphi(f))\le\mathcal{E}(f)$.
\item (Regular property) $\mathcal{F}\cap C_c(X)$ is uniformly dense in $C_c(X)$ and $(\mathcal{E}(\cdot)+\lVert {\cdot}\rVert^p_{L^p(X;m)})^{1/p}$-dense in $\mathcal{F}$.
\item (Strongly local property) For any $f,g\in\mathcal{F}$ with compact support and $g$ constant in an open neighborhood of $\mathrm{supp}(f)$, we have $\mathcal{E}(f+g)=\mathcal{E}(f)+\mathcal{E}(g)$.
\item ($p$-Clarkson's inequality) For any $f,g\in\mathcal{F}$, we have
\begin{equation*}\label{eq_Cla}\tag*{Cla}
\begin{cases}
\mathcal{E}(f+g)+\mathcal{E}(f-g)\ge2 \left(\mathcal{E}(f)^{\frac{1}{p-1}}+\mathcal{E}(g)^{\frac{1}{p-1}}\right)^{p-1}&\text{if }p\in(1,2],\\
\mathcal{E}(f+g)+\mathcal{E}(f-g)\le2 \left(\mathcal{E}(f)^{\frac{1}{p-1}}+\mathcal{E}(g)^{\frac{1}{p-1}}\right)^{p-1}&\text{if }p\in[2,+\infty).\\
\end{cases}
\end{equation*}
\end{enumerate}
We also assume the following condition.
\begin{itemize}
\item ($\mathcal{F}\cap L^\infty(X;m)$ is an algebra) For any $f,g\in\mathcal{F}\cap L^\infty(X;m)$, we have $fg\in\mathcal{F}$ and
\begin{equation*}\label{eq_Alg}\tag*{Alg}
\mathcal{E}(fg)^{{1}/{p}}\le \lVert {f}\rVert_{L^\infty(X;m)}\mathcal{E}(g)^{1/p}+\lVert {g}\rVert_{L^\infty(X;m)}\mathcal{E}(f)^{1/p}.
\end{equation*}
\end{itemize}
For any $\lambda>0$, denote
$$\mathcal{E}_\lambda(\cdot):=\mathcal{E}(\cdot)+\lambda\lVert\cdot \rVert^p_{L^p(X;m)}.$$
The generalized $p$-contraction property introduced in \cite{KS24aV5} implies \ref{eq_Cla} and \ref{eq_Alg} and holds on a large family of metric measure spaces.

By \cite[Theorem 2.4]{Sas26}, a $p$-energy $(\mathcal{E},\mathcal{F})$ corresponds to a (canonical) $p$-energy measure $\Gamma:\mathcal{F}\times\mathcal{B}(X)\to[0,+\infty)$, $(f,A)\mapsto\Gamma(f)(A)$ satisfying the following conditions.
\begin{enumerate}[label=(\alph*),ref=(\alph*)]
\item\label{item_meas1} For any $f\in\mathcal{F}$, $\Gamma(f)(\cdot)$ is a positive Radon measure on $X$ with $\Gamma(f)(X)=\mathcal{E}(f)$.
\item\label{item_meas2} For any $A\in\mathcal{B}(X)$, $\Gamma(\cdot)(A)^{1/p}$ is a seminorm on $\mathcal{F}$.
\item\label{item_meas3} For any $f,g\in\mathcal{F}\cap C_c(X)$, $A\in\mathcal{B}(X)$, if $f-g$ is constant on $A$, then $\Gamma(f)(A)=\Gamma(g)(A)$.
\item\label{item_meas4} ($p$-Clarkson's inequality) For any $f,g\in\mathcal{F}$, for any $A\in\mathcal{B}(X)$, we have
\begin{equation*}
\begin{cases}
\Gamma(f+g)(A)+\Gamma(f-g)(A)\ge2 \left(\Gamma(f)(A)^{\frac{1}{p-1}}+\Gamma(g)(A)^{\frac{1}{p-1}}\right)^{p-1}&\text{if }p\in(1,2],\\
\Gamma(f+g)(A)+\Gamma(f-g)(A)\le2 \left(\Gamma(f)(A)^{\frac{1}{p-1}}+\Gamma(g)(A)^{\frac{1}{p-1}}\right)^{p-1}&\text{if }p\in[2,+\infty).\\
\end{cases}
\end{equation*}
\item\label{item_meas5} (Chain rule) For any $f,g\in\mathcal{F}\cap C_c(X)$, for any piecewise $C^1$ functions $\varphi,\psi:\mathbb{R}\to\mathbb{R}$ with $\varphi(0)=\psi(0)=0$, we have
$$\mathrm{d}\Gamma\left(\varphi(f);\psi(g)\right)=\mathrm{sgn}(\varphi'(f))|\varphi'(f)|^{p-1}\psi'(g)\mathrm{d}\Gamma(f;g),$$
where $\Gamma(f;g)$ is a signed measure given by $\Gamma(f;g)=\frac{1}{p}\frac{\mathrm{d}}{\mathrm{d}t}\Gamma(f+tg)|_{t=0}$.
\item\label{item_meas6} (Leibniz rule) For any $f,g,h\in \mathcal{F}\cap C_c(X)$, we have $\mathrm{d}\Gamma(f;gh)=g \mathrm{d}\Gamma(f;h)+h \mathrm{d}\Gamma(f;g)$.
\end{enumerate}
Using the chain rule, we have the following condition.
\begin{itemize}
\item (Strong sub-additivity) For any $f,g\in\mathcal{F}$, we have $f\vee g, f\wedge g\in\mathcal{F}$ and
\begin{equation*}\label{eq_SubAdd}\tag*{SubAdd}
\mathcal{E}(f\vee g)+\mathcal{E}(f\wedge g)\le\mathcal{E}(f)+\mathcal{E}(g).
\end{equation*}
\end{itemize}

By \ref{eq_Cla} and \cite[Theorem 3.7]{KS24aV5}, for any $u,v\in\mathcal{F}$, the derivative
$$\mathcal{E}(u;v):=\frac{1}{p}\left.\frac{\mathrm{d}}{\mathrm{d}t}\mathcal{E}(u+tv)\right|_{t=0}\in \mathbb{R}$$
exists, the map $\mathcal{E}(u;\cdot):\mathcal{F}\to\mathbb{R}$ is linear, and $\mathcal{E}(u;u)=\mathcal{E}(u)$.


For any bounded open subset $U\subseteq X$, let $\mathcal{F}(U)$ be the closure of $\mathcal{F}\cap C_c(U)$ with respect to $\mathcal E_1^{1/p}$. By \cite[Lemma 6.2]{Yan25c},
$$\mathcal{F}(U)=\left\{u\in \mathcal{F}:\widetilde{u}=0\text{ q.e. on }X\backslash U\right\},$$
where $\widetilde{u}$ is a quasi-continuous modification of $u$; see \cite[Section 8]{Yan25a} for more details.

We say that $u\in\mathcal{F}$ is harmonic in $U$ if $\mathcal{E}(u;v)=0$ for any $v\in\mathcal{F}\cap C_c(U)$. We say that $u$ is superharmonic in $U$ (resp. subharmonic in $U$) if $\mathcal{E}(u;v)\ge0$ (resp. $\mathcal E(u;v)\le0$) for any non-negative $v\in\mathcal{F}\cap C_c(U)$. Since $U$ is bounded, these definitions are unchanged when $\mathcal{F}\cap C_c(U)$ is replaced by $\mathcal{F}(U)$, retaining the non-negativity requirement in the superharmonic and subharmonic cases.

Let $\beta_p>0$. We say that the Poincar\'e inequality \ref{eq_PI} holds if there exist constants $C_{PI}>0$ and $A_{PI}\ge1$ such that for any ball $B(x,r)$ and any $f\in\mathcal{F}$,
\begin{equation*}\label{eq_PI}\tag*{$\text{PI}(\beta_p)$}
\int_{B(x,r)}\lvert f-f_{B(x,r)}\rvert^p\mathrm{d} m\le C_{PI}r^{\beta_p}\int_{{B(x,A_{PI}r)}}\mathrm{d}\Gamma(f).
\end{equation*}

Let $A_1,A_2\in\mathcal{B}(X)$. We define the $\mathcal{E}$-capacity between $A_1$ and $A_2$ as
\begin{align*}
&\mathrm{cap}(A_1,A_2):=\inf\left\{\mathcal{E}(\phi):\phi\in\mathcal{F},
\begin{array}{l}
\phi=1\text{ in an open neighborhood of }A_1,\\
\phi=0\text{ in an open neighborhood of }A_2
\end{array}
\right\},
\end{align*}
where we use the convention that $\inf\emptyset=+\infty$.

We say that the capacity upper bound \ref{eq_ucap} holds if there exist constants $C_{cap}>0$ and $A_{cap}>1$ such that for any ball $B(x,r)$,
\begin{equation*}\label{eq_ucap}\tag*{$\text{cap}(\beta_p)_{\le}$}
\mathrm{cap}\left(B(x,r),X\backslash B(x,A_{cap}r)\right)\le C_{cap} \frac{V(x,r)}{r^{\beta_p}}.
\end{equation*}

\begin{figure}[ht]
  \centering
  \begin{tikzpicture}[scale=1.1]
    \def\pheight{2} 
    \def\xmax{4.5}

    \fill[blue,opacity=0.5]
      (1,\pheight)
      -- (\xmax,\pheight)
      -- (\xmax,{\xmax+\pheight-1})
      -- cycle;

    \draw[blue,opacity=0.5,thick]
      (1,\pheight) -- (\xmax,\pheight);
    \draw[blue,opacity=0.5,thick]
      (1,\pheight) -- (\xmax,{\xmax+\pheight-1});

    \draw[black,dashed]
      (0,\pheight) -- (1,\pheight);
    \draw[black,dashed]
      (0,{\pheight-1}) -- (1,\pheight);
    \draw[black,dashed]
      (1,0) -- (1,\pheight);

\draw[->] (0,0) -- (5.2,0) node[below] {$d_h$};
\draw[->] (0,0) -- (0,6) node[left] {$\beta_p$};

    \draw (1,0.06) -- (1,-0.06) node[below] {$1$};
    \draw (0.06,\pheight) -- (-0.06,\pheight)
      node[left] {$p$};
    \draw (0.06,{\pheight-1}) -- (-0.06,{\pheight-1})
      node[left] {$p-1$};

  \end{tikzpicture}
  \caption{The admissible region of $(d_h,\beta_p)$}
  \label{fig_betapdh}
\end{figure}
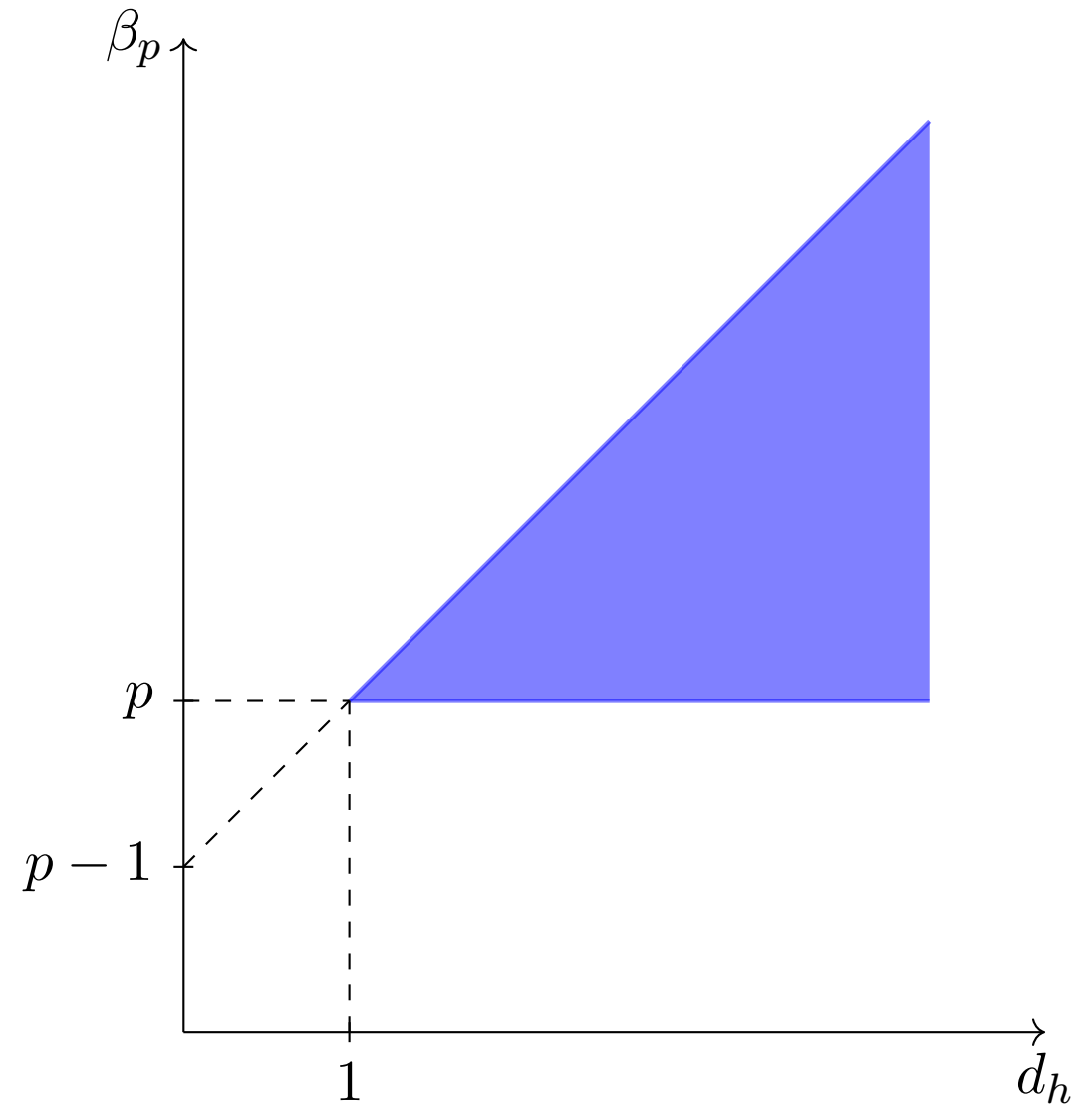

Fix $p>1$, assume \ref{eq_VD}, and let $(\mathcal{E},\mathcal{F})$ be a $p$-energy satisfying \ref{eq_PI} and \ref{eq_ucap}. By \cite[Proposition 2.1 and Remark 2.2]{Yan25a},
\begin{equation}\label{eq_betap_lbd}
\beta_p\ge p.
\end{equation}
By \cite{Eri26}, the substantially stronger cutoff Sobolev inequality $\text{CS}(\beta_p)$ holds. Consequently, when applying results from the literature that assume \ref{eq_VD}, \ref{eq_PI}, and $\text{CS}(\beta_p)$, it suffices to assume \ref{eq_VD}, \ref{eq_PI}, and \ref{eq_ucap}. If \ref{eq_Vdh} also holds, then \cite[Proposition 2.1]{Yan25a} gives
\begin{equation}\label{eq_betap_bdd}
p\le\beta_p\le d_h+p-1,
\end{equation}
see Figure \ref{fig_betapdh}.

In what follows, we use the superscript $(p)$ to indicate explicit dependence on $p$. For example, we write $(\mathcal{E}^{(p)},\mathcal{F}^{(p)})$ for a $p$-energy, \hyperref[eq_PI]{$\text{PI}^{(p)}(\cdot)$} for the Poincar\'e inequality, and $C^{(p)}$ and $A^{(p)}$ for constants depending on $p$.

The first main result of this paper is the following dichotomy.

\begin{theorem}\label{thm_dichotomy1}
Assume \ref{eq_VD}. Let $I\subseteq(1,+\infty)$ be an open interval. Assume that, for any $p\in I$, there exists a $p$-energy $(\mathcal{E}^{(p)},\mathcal{F}^{(p)})$ satisfying \hyperref[eq_PI]{$\text{PI}^{(p)}(\beta_p)$} and \hyperref[eq_ucap]{$\text{cap}^{(p)}(\beta_p)_\le$}, where $\beta_p\ge p$ by (\ref{eq_betap_lbd}). Then exactly one of the following alternatives holds:
\begin{enumerate}[label=(\roman*),ref=(\roman*)]
\item\label{item_dich1} $\beta_p=p$ for any $p\in I$;
\item\label{item_dich2} the function $p\mapsto \frac{\beta_p}{p}$ takes values in $(1,+\infty)$ and is \emph{strictly} decreasing on $I$.
\end{enumerate}
\end{theorem}

\begin{remark}
Roughly speaking, $\beta_p$ can be interpreted as the scaling order of the $p$-Laplacian. In Euclidean space and other settings with a classical gradient structure, we typically have $\beta_p=p$, so that $\frac{\beta_p}{p}=1$ corresponds to the order of a first derivative. On fractals, we typically have $\beta_p>p$, and the ratio $\frac{\beta_p}{p}$ may instead be viewed heuristically as the ``effective" order of an analogue of the gradient. Its deviation from 1 thus reflects a departure from classical first-order scaling. In this sense, the above result means that the scaling of the $p$-energy becomes ``progressively" closer to classical first-order scaling as $p$ increases.
\end{remark}

\begin{remark}
By \cite[Theorem 2.1, Lemmas 5.1 and 5.2]{Yan25e}, the dichotomy holds with \ref{item_dich2} replaced by the following weaker alternative:
\begin{enumerate}[label=(\roman*${}'$),ref=(\roman*${}'$)]
\setcounter{enumi}{1}
\item\label{item_dich2p} the function $p\mapsto\frac{\beta_p}{p}$ takes values in $(1,+\infty)$ and is decreasing on $I$.
\end{enumerate}
It therefore remains only to prove that the monotonicity in \ref{item_dich2p} is strict.
\end{remark}

\begin{remark}\label{rmk_Fen}
We sketch how the strict monotonicity result relates to the unboundedness of the Riesz transform for any $p>2$ on the fractal-like spaces considered in \cite{Fen26}. The following discussion is heuristic and does not constitute a formal proof.

Assume that on the corresponding fractals, for any $p\in(1,2]$, there exists a $p$-energy satisfying \hyperref[eq_PI]{$\text{PI}^{(p)}(\beta_p)$} and \hyperref[eq_ucap]{$\text{cap}^{(p)}(\beta_p)_\le$}. Then suitable variants of these two inequalities also hold on the fractal-like graphs and manifolds.

In the graph or manifold setting considered there, a suitable gradient operator $\nabla$ is well-defined. We use the notation of \cite{Fen26} and refer the reader to that paper for further details. By duality, to establish the unboundedness of the Riesz transform for any $p>2$, it suffices to prove the failure of the reverse Riesz inequality $(\text{RR}_{p,\frac{1}{2}})$ for any $p\in(1,2)$. Under the sub-Gaussian upper heat kernel estimate $(\text{UE}_\beta)$ with $\beta>2$, \cite[Proposition 4.6 and Remark 4.7]{Fen26} show that $(\text{RR}_{p,\frac{1}{2}})$ implies $(\text{P}_{p,\frac{\beta}{2}})$, that is, for any ball $B=B(x,r)$ and any $f$ supported in $B$,
\begin{equation}\label{eq_Fen1}
\lVert f\rVert_{L^p(B)}\lesssim r^{\frac{\beta}{2}}\lVert \nabla f\rVert_{L^p(B)}.
\end{equation}
Here, $\beta$ in $(\text{UE}_\beta)$ corresponds to the $2$-walk dimension $\beta_2$ in our setting. The transferred version of \hyperref[eq_ucap]{$\text{cap}^{(p)}(\beta_p)_\le$} ensures that for any ball $B=B(x,r)$ with $r>1$, there exists a function $f$ supported in $B$ such that $0\le f\le 1$ in $B$, $f=1$ in $\frac{1}{2}B$, and
\begin{equation}\label{eq_Fen2}
\lVert \nabla f\rVert_{L^p(B)}^p \lesssim \frac{m(B)}{r^{\beta_p}}.
\end{equation}
Moreover,
$$\lVert f\rVert_{L^p(B)}^p\ge m\left(\frac{1}{2}B\right)\asymp m(B),$$
where the comparability follows from the volume doubling condition \ref{eq_VD} on the fractal-like space. Combining this with (\ref{eq_Fen1}) and (\ref{eq_Fen2}), we obtain
\begin{equation}\label{eq_Fen3}
r^{\frac{\beta_p}{p}}\lesssim r^{\frac{\beta_2}{2}}\text{ for any $r>1$}.
\end{equation}
Since $\beta_2>2$, Theorem \ref{thm_dichotomy1} gives $\frac{\beta_p}{p}>\frac{\beta_2}{2}$ for any $p\in(1,2)$. Letting $r\to+\infty$ in (\ref{eq_Fen3}), we obtain a contradiction.
\end{remark}

The second main result of this paper is as follows.

\begin{theorem}\label{thm_gammap}
Assume \ref{eq_Vdh}. Let $I\subseteq(1,+\infty)$ be an open interval. Assume that, for any $p\in I$, there exists a $p$-energy $(\mathcal{E}^{(p)},\mathcal{F}^{(p)})$ satisfying \hyperref[eq_PI]{$\text{PI}^{(p)}(\beta_p)$} and \hyperref[eq_ucap]{$\text{cap}^{(p)}(\beta_p)_\le$}, where $\beta_p\le d_h+p-1$ by (\ref{eq_betap_bdd}). Then:
\begin{enumerate}[label=(\arabic*),ref=(\arabic*)]
\item\label{item_gammap1} The function $p\mapsto \gamma_p:=\frac{\beta_p-d_h}{p-1}$ is increasing on $I$.
\item\label{item_gammap2} If $\gamma_p=1$ for some $p\in I$, then $\gamma_q=1$ for any $q\in I$.
\item\label{item_gammap3} Let $p,q\in I$ satisfy $p<q$. If $\gamma_p<1$, or equivalently, $\gamma_q<1$, then $\gamma_p<\gamma_q$.
\end{enumerate}
Consequently, exactly one of the following alternatives holds:
\begin{enumerate}[label=(\alph*)]
\item $\beta_p=d_h+p-1$ for any $p\in I$;
\item the function $p\mapsto \frac{\beta_p-d_h}{p-1}$ takes values in $(-\infty,1)$ and is \emph{strictly} increasing on $I$.
\end{enumerate}
\end{theorem}

\begin{remark}
Roughly speaking, for $p$ sufficiently large that $\beta_p>d_h$ and $\gamma_p\in(0,1]$, the exponent $\gamma_p$ describes the H\"older regularity of $p$-harmonic functions (see \cite[Theorem 6.31]{KS24aV5}), which may be viewed heuristically as among the ``most" regular functions in the domain of the $p$-energy. The case $\gamma_p=1$ occurs in spaces with some kind of tree structure, such as $\mathbb{R}$ and the Vicsek set (see \cite{BC23}), where $p$-harmonic functions are indeed Lipschitz. If $\gamma_p<1$ for some, or equivalently any, admissible $p$, the above result and (\ref{eq_betap_bdd}) imply that $\gamma_p$ is strictly increasing and converges to $1$ as $p\to+\infty$. Hence, $p$-harmonic functions become ``strictly smoother" as $p$ increases, in the sense that their H\"older regularity exponents strictly increase towards the Lipschitz exponent $1$.
\end{remark}

As a direct corollary, we obtain the following strict monotonicity result, which covers both the concrete fractals discussed in the Introduction and the Sierpi\'nski carpet.

\begin{corollary}
Let $K$ be the Sierpi\'nski gasket, the Sierpi\'nski carpet, or the golden ratio Sierpi\'nski gasket, and let $d_h$ be its Hausdorff dimension. For any $p\in(1,+\infty)$, let $(\mathcal{E}^{(p)},\mathcal{F}^{(p)})$ be the $p$-energy constructed in \cite{HPS04}, \cite{MS25}, or \cite{ChenQiu26}, respectively, with $p$-walk dimension $\beta_p$. Then
\begin{itemize}
\item the function $p\mapsto \frac{\beta_p}{p}$ takes values in $(1,+\infty)$ and is strictly decreasing on $(1,+\infty)$;
\item the function $p\mapsto \frac{\beta_p-d_h}{p-1}$ takes values in $(-\infty,1)$ and is strictly increasing on $(1,+\infty)$.
\end{itemize}
\end{corollary}

\begin{proof}
On the Sierpi\'nski gasket, the resistance estimate \hyperref[eq_res]{$\text{R}^{(p)}(\beta_p)$}, defined below, holds by \cite[Theorem 4.6]{Kig23} and \cite[Theorem 8.30]{KS24aV5}. On the golden ratio Sierpi\'nski gasket, the same estimate holds by \cite[Theorem 2]{ChenQiu26}. By \cite[Theorem 2.3]{Yan25c}, this resistance estimate is equivalent to the conjunction of \hyperref[eq_PI]{$\text{PI}^{(p)}(\beta_p)$} and \hyperref[eq_ucap]{$\text{cap}^{(p)}(\beta_p)_\le$}. Hence, both inequalities hold on these two fractals. On the Sierpi\'nski carpet, these two inequalities hold by \cite[Theorem 8.21 and Proposition 6.21]{MS25}. By the well-known fact that $2<\beta_2<d_h+1$ on these fractals, the results follow directly from Theorems \ref{thm_dichotomy1} and \ref{thm_gammap}.
\end{proof}

This paper is organized as follows. In Section \ref{sec_pre}, we collect some preliminary results for later use. In Section \ref{sec_dichotomy1}, we give the proof of Theorem \ref{thm_dichotomy1}. In Section \ref{sec_gammap1}, we give the proof of Theorem \ref{thm_gammap}\ref{item_gammap1}. In Section \ref{sec_gammap2}, we give the proof of Theorem \ref{thm_gammap}\ref{item_gammap2}. In Section \ref{sec_gammap3}, we give the proof of Theorem \ref{thm_gammap}\ref{item_gammap3}.

\section{Preliminaries}\label{sec_pre}

In this section, we collect some preliminary results for later use.

For any $p>1$, $\beta>0$, and $r>0$, let $U\subseteq X$ be open and let $u$ be a measurable function. We define
$$E^{(p)}_{\beta,U}(u,r):=\frac{1}{r^\beta}\int_U \left(\frac{1}{V(x,r)}\int_{B(x,r)}\lvert u(x)-u(y)\rvert^p m(\mathrm{d}y)\right)m(\mathrm{d}x),$$
and
$$U(r):=\bigcup_{x\in U}B(x,r).$$

For any $\alpha>0$, define the Besov space
$$B^{p,\alpha}(X):=\left\{u\in L^p(X;m):\sup_{r>0} E^{(p)}_{p\alpha,X}(u,r)<+\infty\right\}.$$
Obviously, $B^{p,\alpha}(X)$ is decreasing in $\alpha$ and may become trivial if $\alpha$ is too large. We define the following critical exponent
\begin{align*}
&\alpha_p(X):=\sup\left\{\alpha>0:B^{p,\alpha}(X)\text{ contains non-constant functions}\right\}\le+\infty.
\end{align*}
The value of $\alpha_p(X)$ depends \emph{only} on the metric measure space $(X,d,m)$. The following lemma records some basic properties of $\alpha_p(X)$.

\begin{lemma}[{\cite[Theorem 4.1]{Bau24}}]\label{lem_Bau}
Assume \ref{eq_VD}.
\begin{enumerate}[label=(\roman*)]
\item For any $p\in(1,+\infty)$, we have $\alpha_p(X)\ge1$.
\item The function $p\mapsto p\alpha_p(X)$ is monotone increasing for $p\in(1,+\infty)$.
\item The function $p\mapsto \alpha_p(X)$ is monotone decreasing for $p\in(1,+\infty)$.
\end{enumerate}
Hence
\begin{enumerate}[label=(\alph*)]
\item For $p\in(1,+\infty)$, the functions $p\mapsto p\alpha_p(X)$ and $p\mapsto \alpha_p(X)$ are continuous.
\item If $\alpha_p(X)<+\infty$ for some $p\in(1,+\infty)$, then $\alpha_p(X)<+\infty$ for any $p\in(1,+\infty)$.
\end{enumerate}
\end{lemma}

The value of $\alpha_p(X)$ can be determined provided that there exists a $p$-energy satisfying certain functional inequalities.

\begin{lemma}[{\cite[Theorem 4.6]{Shi24a}}]\label{lem_Shi1}
Assume \ref{eq_VD}. Let $(\mathcal{E}^{(p)},\mathcal{F}^{(p)})$ be a $p$-energy satisfying \hyperref[eq_PI]{$\text{PI}^{(p)}(\beta_p)$} and \hyperref[eq_ucap]{$\text{cap}^{(p)}(\beta_p)_\le$}. Then $\alpha_p(X)=\frac{\beta_p}{p}$.
\end{lemma}

We have the following results concerning the quantity $E^{(p)}_{\beta,U}$.

\begin{lemma}[{\cite[Proposition 3.5 and Equation (1.18)]{Shi24a}}]\label{lem_Shi2}
Assume \ref{eq_VD}. Let $(\mathcal{E}^{(p)},\mathcal{F}^{(p)})$ be a $p$-energy satisfying \hyperref[eq_PI]{$\text{PI}^{(p)}(\beta_p)$}. Then there exists a constant $C^{(p)}>0$ such that for any open subset $U\subseteq X$, $u\in \mathcal{F}^{(p)}$, and $r>0$,
$$E^{(p)}_{\beta_p,U}(u,r)\le C^{(p)}\Gamma^{(p)}(u)\left(U(2A_{PI}^{(p)}r)\right),$$
where $A^{(p)}_{PI}\ge1$ is the constant in \hyperref[eq_PI]{$\text{PI}^{(p)}(\beta_p)$}.
\end{lemma}

\begin{lemma}[{\cite[Corollaries 1.14 and 1.15, and Equation (1.18)]{Shi24a}}]\label{lem_Shi3}
Assume \ref{eq_VD}. Let $(\mathcal{E}^{(p)},\mathcal{F}^{(p)})$ be a $p$-energy satisfying \hyperref[eq_PI]{$\text{PI}^{(p)}(\beta_p)$} and \hyperref[eq_ucap]{$\text{cap}^{(p)}(\beta_p)_\le$}. Then
\begin{align*}
&\mathcal{F}^{(p)}=\left\{u\in L^p(X;m):\liminf_{r\downarrow0}E^{(p)}_{\beta_p,X}(u,r)<+\infty\right\}\\
&=\left\{u\in L^p(X;m):\limsup_{r\downarrow0}E^{(p)}_{\beta_p,X}(u,r)<+\infty\right\}=B^{p,\frac{\beta_p}{p}}(X).
\end{align*}
Moreover, there exists a constant $C^{(p)}\ge1$ such that for any $u\in \mathcal{F}^{(p)}$,
$$\frac{1}{C^{(p)}}\mathcal{E}^{(p)}(u)\le \liminf_{r\downarrow0}E^{(p)}_{\beta_p,X}(u,r)\le \sup_{r>0}E^{(p)}_{\beta_p,X}(u,r)\le C^{(p)}\mathcal{E}^{(p)}(u).$$
\end{lemma}

\begin{lemma}[{\cite[Theorem 3.12 and Equation (1.18)]{Shi24a}}]\label{lem_Shi4}
Assume \ref{eq_VD}. Let $(\mathcal{E}^{(p)},\mathcal{F}^{(p)})$ be a $p$-energy satisfying \hyperref[eq_PI]{$\text{PI}^{(p)}(\beta_p)$} and \hyperref[eq_ucap]{$\text{cap}^{(p)}(\beta_p)_\le$}. Then there exists a constant $C^{(p)}>0$ such that for any open subset $U\subseteq X$ and any $u\in \mathcal{F}^{(p)}$,
$$\Gamma^{(p)}(u)(U)\le C^{(p)}\liminf_{r\downarrow0}E^{(p)}_{\beta_p,U}(u,r).$$
\end{lemma}

Let $A_1,A_2\in\mathcal{B}(X)$. We define the $\mathcal{E}^{(p)}$-effective resistance between $A_1$ and $A_2$ as
$$R^{(p)}(A_1,A_2):=\frac{1}{\mathrm{cap}^{(p)}(A_1,A_2)},$$
where we use the convention that $0^{-1}=+\infty$  and $(+\infty)^{-1}=0$. If $A_1\subseteq A_2$ and $B_1\subseteq B_2$, then
$$R^{(p)}(A_1,B_1)\ge R^{(p)}(A_2,B_2).$$
For any $x,y\in X$ with $x\ne y$, let $R^{(p)}(x,x)=R^{(p)}(\{x\},\{x\})=0$ and
$$R^{(p)}(x,y)=R^{(p)}(\{x\},\{y\}):=\lim_{r\downarrow0}R^{(p)}\left(B(x,r),B(y,r)\right).$$

Let $\beta_p>0$. We say that the resistance estimate \ref{eq_res} holds if there exists a constant $C_R\ge1$ such that for any $x,y\in X$ with $x\ne y$,
\begin{equation*}\label{eq_res}\tag*{$\text{R}^{(p)}(\beta_p)$}
\frac{1}{C_R}\frac{d(x,y)^{\beta_p}}{m(B(x,d(x,y)))}\le R^{(p)}(x,y)\le C_R\frac{d(x,y)^{\beta_p}}{m(B(x,d(x,y)))}.
\end{equation*}

\begin{lemma}\label{lem_rec}
Assume \ref{eq_Vdh} and $\beta_p>d_h$. Let $(\mathcal{E}^{(p)},\mathcal{F}^{(p)})$ be a $p$-energy. Then
\begin{equation}\label{eq_equiv_res}
\text{\hyperref[eq_PI]{$\text{PI}^{(p)}(\beta_p)$}}+\text{\hyperref[eq_ucap]{$\text{cap}^{(p)}(\beta_p)_\le$}}\Longleftrightarrow \text{\hyperref[eq_res]{$\text{R}^{(p)}(\beta_p)$}}.
\end{equation}
Moreover, whenever these equivalent conditions hold, the following assertions hold:
\begin{enumerate}[label=(\arabic*),ref=(\arabic*)]
\item\label{item_rho} $\mathcal{F}^{(p)}\subseteq C(X)$. For any $x,y\in X$ with $x\ne y$,
$$R^{(p)}(x,y)=\sup\left\{\frac{|u(x)-u(y)|^p}{\mathcal{E}^{(p)}(u)}:u\in\mathcal{F}^{(p)}\cap C_c(X),\ \mathcal{E}^{(p)}(u)>0\right\},$$
which implies
$$|u(x)-u(y)|^p\le R^{(p)}(x,y)\mathcal{E}^{(p)}(u)\text{ for any }u\in\mathcal{F}^{(p)},\ x,y\in X.$$
Moreover,
\begin{equation}\label{eq_rhop}
\rho^{(p)}:=(R^{(p)})^{\frac{1}{p-1}}:X\times X\to [0,+\infty),\qquad (x,y)\mapsto (R^{(p)}(x,y))^{\frac{1}{p-1}}
\end{equation}
is a metric on $X$, and there exists a constant $L\ge1$ such that
\begin{equation}\label{eq_rhodc}
\frac{1}{L}d(x,y)^{\frac{\beta_p-d_h}{p-1}}\le\rho^{(p)}(x,y)\le{L}d(x,y)^{\frac{\beta_p-d_h}{p-1}}\text{ for any }x,y\in X.
\end{equation}
We write $B_{\rho^{(p)}}(x,r)$ for the ball centered at $x$ with radius $r$ with respect to $\rho^{(p)}$, and denote by $\mathcal{H}^1_{\rho^{(p)}}$ the $1$-dimensional Hausdorff measure on $(X,\rho^{(p)})$.
\item\label{item_pot} There exist constants $C_1,C_2>0$ such that the following holds. For any compact set $K$ and any bounded open set $U$ satisfying $K\subseteq U$, let $u\in \mathcal{F}^{(p)}\subseteq C(X)$ be the $\mathcal{E}^{(p)}$-capacitary potential for $(K,X\backslash U)$; that is, $0\le u\le 1$ on $X$, $u=1$ on $K$, $u=0$ on $X\backslash U$, $u$ is harmonic in $U\backslash K$, and $\mathcal{E}^{(p)}(u)=\mathrm{cap}^{(p)}(K,X\backslash U)$. Then
\begin{align}
&\lvert u(x)-u(y)\rvert\le C_1 \mathcal{E}^{(p)}(u)^{\frac{1}{p-1}}{d(x,y)^{\frac{\beta_p-d_h}{p-1}}}\label{eq_uholder}\\
&\le C_2 \mathcal{E}^{(p)}(u)^{\frac{1}{p-1}}{\rho^{(p)}(x,y)}\text{ for any }x,y\in X.\label{eq_ulip}
\end{align}
Moreover, the pushforward $u_*\Gamma^{(p)}(u)$ of $\Gamma^{(p)}(u)$ under $u$ is a constant multiple of the one-dimensional Lebesgue measure $\mathcal{L}^1$ on $[0,1]$. More precisely,
\begin{equation}\label{eq_upush}
u_{*}\Gamma^{(p)}(u)=\mathcal{E}^{(p)}(u)\mathcal{L}^1|_{[0,1]}.
\end{equation}
Consequently, for any $x\in X$ and any $r>0$,
\begin{equation}\label{eq_mu_upper}
\Gamma^{(p)}(u)(B_{\rho^{(p)}}(x,r))\le  2C_2\mathcal{E}^{(p)}(u)^{\frac{p}{p-1}}r,
\end{equation}
and
\begin{equation}\label{eq_mu_upperA}
\Gamma^{(p)}(u)\le 2C_2\mathcal{E}^{(p)}(u)^{\frac{p}{p-1}} \mathcal{H}^1_{\rho^{(p)}}.
\end{equation}
\end{enumerate}
\end{lemma}

\begin{proof}
The equivalence (\ref{eq_equiv_res}) follows from \cite[Theorem 2.3]{Yan25c}.

\ref{item_rho} The function $\rho^{(p)}=(R^{(p)})^{\frac{1}{p-1}}$ is a metric by the argument of \cite[Corollary 6.32]{KS24aV5}. Indeed, that argument remains valid under $(\mathrm{RF1})_p$--$(\mathrm{RF4})_p$, the Markovian property, $p$-Clarkson's inequality, and strong sub-additivity; the full generalized $p$-contraction property is not required. These properties also hold for $(\mathcal{E}^{(p)},\mathcal{F}_e^{(p)})$, where $\mathcal{F}_e^{(p)}$ will be given in (\ref{eq_Fpe}). The remaining results follow from \cite[Lemma 5.2]{Yan25c}, \hyperref[eq_res]{$\text{R}^{(p)}(\beta_p)$}, and the same arguments as those used in \cite[Lemma 8.6 and Corollary 8.9]{GY26}.

\ref{item_pot} The existence of such an $\mathcal{E}^{(p)}$-capacitary potential $u$ follows from \cite[Lemma 3.8]{Yan25d}. We claim that for any $a,b\in[0,1]$ with $a<b$,
\begin{equation}\label{eq_Gammaab}
\Gamma^{(p)}(u)\left(\left\{x\in X:a<u(x)<b\right\}\right)=(b-a)\mathcal{E}^{(p)}(u).
\end{equation}
Indeed, since $(u-a)_+\wedge (b-a)-(b-a)u\in \mathcal{F}^{(p)}(U\backslash K)$ and $u$ is harmonic in $U\backslash K$, by the strongly local property,
\begin{align*}
&0=\mathcal{E}^{(p)}(u;(u-a)_+\wedge (b-a))-(b-a)\mathcal{E}^{(p)}(u)\\
&=\Gamma^{(p)}(u)\left(\{x\in X:a<u(x)<b\}\right)-(b-a)\mathcal{E}^{(p)}(u).
\end{align*}

For any $x,y\in X$ with $x\ne y$, without loss of generality, we may assume that $u(x)>u(y)$, let $v=(u-u(y))_+\wedge (u(x)-u(y))\in \mathcal{F}^{(p)}\subseteq C(X)$. By (\ref{eq_Gammaab}),
$$R^{(p)}(x,y)\ge \frac{\lvert v(x)-v(y)\rvert^p}{\mathcal{E}^{(p)}(v)}=\frac{\lvert u(x)-u(y)\rvert^p}{\Gamma^{(p)}(u)\left(\{z\in X:u(y)<u(z)<u(x)\}\right)}=\frac{\lvert u(x)-u(y)\rvert^{p-1}}{\mathcal{E}^{(p)}(u)}.$$
Combining this with \hyperref[eq_res]{$\text{R}^{(p)}(\beta_p)$},
$$\lvert u(x)-u(y)\rvert^{p-1}\lesssim \mathcal{E}^{(p)}(u)d(x,y)^{\beta_p-d_h},$$
which gives (\ref{eq_uholder}). Together with (\ref{eq_rhodc}), we obtain (\ref{eq_ulip}).

To prove (\ref{eq_upush}), take any $\phi\in C^1([0,1])$ satisfying $\phi(0)=\phi(1)=0$. For any $\varepsilon\in \mathbb{R}$, let $\psi(t)=t+\varepsilon \phi(t)$. Then $\psi\circ u\in \mathcal{F}^{(p)}$, with $\psi\circ u=1$ on $K$ and $\psi\circ u=0$ on $X\backslash U$. By the chain rule,
\begin{align*}
&\int_{[0,1]}\mathrm{d}\left(u_{*}\Gamma^{(p)}(u)\right)=\mathcal{E}^{(p)}(u)=\mathrm{cap}^{(p)}(K,X\backslash U)\le \mathcal{E}^{(p)}(\psi\circ u)\\
&=\int_X \lvert\psi'(u)\rvert^p \mathrm{d}\Gamma^{(p)}(u)=\int_{[0,1]}\lvert 1+\varepsilon \phi'(t)\rvert^p \left(u_{*}\Gamma^{(p)}(u)\right)(\mathrm{d}t).
\end{align*}
Since $\varepsilon\in \mathbb{R}$ is arbitrary and $\phi'$ is bounded on $[0,1]$, differentiating with respect to $\varepsilon$ and then letting $\varepsilon=0$, we have
$$\int_{[0,1]}\phi' \mathrm{d}\left(u_{*}\Gamma^{(p)}(u)\right)=0.$$
By Lemma \ref{lem_Lebesgue} below,
$$u_{*}\Gamma^{(p)}(u)=\left(\left(u_{*}\Gamma^{(p)}(u)\right)([0,1])\right) \mathcal{L}^1|_{[0,1]}=\mathcal{E}^{(p)}(u)\mathcal{L}^1|_{[0,1]}.$$
For any ball $B_{\rho^{(p)}}(x,r)$, by (\ref{eq_ulip}) and (\ref{eq_upush}),
\begin{align*}
&\Gamma^{(p)}(u)(B_{\rho^{(p)}}(x,r))\le \Gamma^{(p)}(u) \left(\left\{y\in X:\lvert u(x)-u(y)\rvert< C_2 \mathcal{E}^{(p)}(u)^{\frac{1}{p-1}}r\right\}\right)\\
&=\mathcal{E}^{(p)}(u)\mathcal{L}^1 \left(\left\{t\in[0,1]:\lvert t-u(x)\rvert<C_2 \mathcal{E}^{(p)}(u)^{\frac{1}{p-1}}r\right\}\right)\le 2C_2\mathcal{E}^{(p)}(u)^{\frac{p}{p-1}}r,
\end{align*}
which is (\ref{eq_mu_upper}). By Frostman's lemma or the mass distribution principle (see \cite[Theorem 8.17]{Mat95book}, \cite[p. 67]{Fal14book}, or \cite[Theorem 2.4]{BDCFMT13}), we obtain (\ref{eq_mu_upperA}).
\end{proof}

\begin{lemma}\label{lem_Lebesgue}
Let $\lambda$ be a non-negative finite Borel measure on $[0,1]$. Assume that for any $\phi\in C^1([0,1])$ with $\phi(0)=\phi(1)=0$, we have $\int_{[0,1]}\phi'\mathrm{d}\lambda=0$. Then $\lambda=\lambda([0,1]) \mathcal{L}^1|_{[0,1]}$.
\end{lemma}

\begin{proof}
We only need to show that for any $f\in C([0,1])$,
$$\int_{[0,1]}f \mathrm{d}\lambda=\lambda([0,1])\int_{[0,1]}f \mathrm{d}\mathcal{L}^1.$$
Indeed, let $g(t)=f(t)-\int_{0}^1 f(s)\mathrm{d}s$, then $g\in C([0,1])$ and $\int_0^1 g(t)\mathrm{d}t=0$. Let $\phi(t)=\int_0^t g(s)\mathrm{d}s$, then $\phi\in C^1([0,1])$, $\phi(0)=\phi(1)=0$, and $\phi'=g$. By assumption,
$$0=\int_{[0,1]}\phi' \mathrm{d}\lambda=\int_{[0,1]}g \mathrm{d}\lambda=\int_{[0,1]}f \mathrm{d}\lambda-\lambda([0,1])\int_{[0,1]}f \mathrm{d}\mathcal{L}^1.$$
\end{proof}

\section{Proof of Theorem \ref{thm_dichotomy1}}\label{sec_dichotomy1}

Suppose that there exist $p,q\in I$ with $p<q$ such that $\frac{\beta_p}{p}=\frac{\beta_q}{q}>1$. Let $U\subseteq X$ be open and let $u$ be a measurable function. For any $r>0$, by H\"older's inequality,
\begin{align}
&E^{(p)}_{\beta_p,U}(u,r)\le \frac{1}{r^{\beta_p}}\int_U \left(\frac{1}{V(x,r)}\int_{B(x,r)}\lvert u(x)-u(y)\rvert^q m(\mathrm{d}y)\right)^{\frac{p}{q}}m(\mathrm{d}x)\nonumber\\
&\le \frac{1}{r^{\beta_p}} \left[\int_U \left(\frac{1}{V(x,r)}\int_{B(x,r)}\lvert u(x)-u(y)\rvert^q m(\mathrm{d}y)\right)m(\mathrm{d}x)\right]^{\frac{p}{q}}m(U)^{1-\frac{p}{q}}\nonumber\\
&=E^{(q)}_{\beta_q,U}(u,r)^{\frac{p}{q}}m(U)^{1-\frac{p}{q}},\label{eq_EUpq}
\end{align}
where the last equality follows from the fact that $\frac{\beta_p}{p}=\frac{\beta_q}{q}$. By \ref{eq_VD} and \hyperref[eq_PI]{$\text{PI}^{(q)}(\beta_q)$}, we may apply Lemma \ref{lem_Shi2} to obtain, for any $u\in\mathcal{F}^{(q)}$,
\begin{equation}\label{eq_EqGamma}
E^{(q)}_{\beta_q,U}(u,r)\le C^{(q)}\Gamma^{(q)}(u)(U(2Ar)),
\end{equation}
where $A=A^{(q)}_{PI}$.

Since $(\mathcal{E}^{(q)},\mathcal{F}^{(q)})$ is regular, we may choose a non-zero function $u\in\mathcal{F}^{(q)}\cap C_c(X)$. Take $B=B(x_0,R)\supseteq \mathrm{supp}(u)$. For any $r\in(0,R)$, by (\ref{eq_EUpq}),
$$E^{(p)}_{\beta_p,X}(u,r)=E^{(p)}_{\beta_p,2B}(u,r)\le E^{(q)}_{\beta_q,2B}(u,r)^{\frac{p}{q}}m(2B)^{1-\frac{p}{q}}=E^{(q)}_{\beta_q,X}(u,r)^{\frac{p}{q}}m(2B)^{1-\frac{p}{q}}.$$
Since $u\in \mathcal{F}^{(q)}$, by \ref{eq_VD}, \hyperref[eq_PI]{$\text{PI}^{(q)}(\beta_q)$}, and \hyperref[eq_ucap]{$\text{cap}^{(q)}(\beta_q)_\le$}, Lemma \ref{lem_Shi3} gives $\liminf_{r\downarrow0}E^{(q)}_{\beta_q,X}(u,r)<+\infty$, which implies $\liminf_{r\downarrow0}E^{(p)}_{\beta_p,X}(u,r)<+\infty$. Since $u\in C_c(X)\subseteq L^p(X;m)$, by \ref{eq_VD}, \hyperref[eq_PI]{$\text{PI}^{(p)}(\beta_p)$}, and \hyperref[eq_ucap]{$\text{cap}^{(p)}(\beta_p)_\le$}, Lemma \ref{lem_Shi3} gives $u\in \mathcal{F}^{(p)}$. Let
$$\mu:=\Gamma^{(q)}(u),$$
then $\mu$ is a finite positive Radon measure. By the strongly local property, $\mathrm{supp}(\mu)\subseteq \mathrm{supp}(u)$.  For any $r\in(0,\frac{1}{8A}R)$, let $W$ be an $r$-net. By \ref{eq_VD}, there exists an integer $N\ge1$ depending only on $C_{VD}$ and $A$ such that
\begin{equation}\label{eq_bdd_over}
\sum_{z\in W}1_{B(z,4Ar)}\le N.
\end{equation}
By (\ref{eq_EUpq}) and (\ref{eq_EqGamma}), for any $z\in W$,
$$E^{(p)}_{\beta_p,B(z,r)}(u,r)\le C \mu(B(z,4Ar))^{\frac{p}{q}}m(B(z,r))^{1-\frac{p}{q}},$$
where $C=(C^{(q)})^{\frac{p}{q}}$, hence
\begin{equation}\label{eq_EpX}
E^{(p)}_{\beta_p,X}(u,r)\le \sum_{z\in W\cap (2B)}E^{(p)}_{\beta_p,B(z,r)}(u,r)\le C\sum_{z\in W\cap (2B)}\mu(B(z,4Ar))^{\frac{p}{q}}m(B(z,r))^{1-\frac{p}{q}}.
\end{equation}
Since $\beta_q>q$, by \cite[Theorem 2.1]{Yan25c}, $\mu\perp m$. For any $\varepsilon>0$, by the regularity of $\mu$ and $m$, there exist a compact set $K$ and an open set $G$ with $K\subseteq G\subseteq B$ such that
\begin{equation}\label{eq_mumKG}
\mu(X\backslash K)<\varepsilon,\qquad m(G)<\varepsilon.
\end{equation}
Let $r_0=\mathrm{dist}(K,X\backslash G)\in(0,+\infty)$. For any $r<\frac{1}{8A}r_0$, if $B(z,4Ar)\cap K\ne\emptyset$, then $B(z,4Ar)\subseteq G$. Write $W\cap (2B)=W_1\sqcup W_2$, where
\begin{align*}
W_1&:=\left\{z\in W\cap (2B):B(z,4Ar)\cap K\ne\emptyset\right\},\\
W_2&:=\left\{z\in W\cap (2B):B(z,4Ar)\cap K=\emptyset\right\},
\end{align*}
then
$$\bigcup_{z\in W_1}B(z,4Ar)\subseteq G,\qquad\bigcup_{z\in W_2}B(z,4Ar)\subseteq X\backslash K.$$
By (\ref{eq_EpX}), H\"older's inequality, (\ref{eq_bdd_over}), and (\ref{eq_mumKG}),
\begin{align*}
&E^{(p)}_{\beta_p,X}(u,r)\le C\sum_{k=1}^2\sum_{z\in W_k}\mu(B(z,4Ar))^{\frac{p}{q}}m(B(z,r))^{1-\frac{p}{q}}\\
&\le C\sum_{k=1}^2\left(\sum_{z\in W_k}\mu(B(z,4Ar))\right)^{\frac{p}{q}}\left(\sum_{z\in W_k}m(B(z,r))\right)^{1-\frac{p}{q}}\\
&\le CN\sum_{k=1}^2\mu \left(\bigcup_{z\in W_k}B(z,4Ar)\right)^{\frac{p}{q}}m \left(\bigcup_{z\in W_k}B(z,r)\right)^{1-\frac{p}{q}}\\
&\le CN \left(\mu(X)^{\frac{p}{q}}m(G)^{1-\frac{p}{q}}+\mu(X\backslash K)^{\frac{p}{q}}m(4B)^{1-\frac{p}{q}}\right)\\
&<CN \left(\mu(X)^{\frac{p}{q}}\varepsilon^{1-\frac{p}{q}}+\varepsilon^{\frac{p}{q}}m(4B)^{1-\frac{p}{q}}\right).
\end{align*}
Hence $\lim_{r\downarrow0}E^{(p)}_{\beta_p,X}(u,r)=0$. By \ref{eq_VD}, \hyperref[eq_PI]{$\text{PI}^{(p)}(\beta_p)$}, and \hyperref[eq_ucap]{$\text{cap}^{(p)}(\beta_p)_\le$}, Lemma \ref{lem_Shi3} again gives $\mathcal{E}^{(p)}(u)=0$. By \hyperref[eq_PI]{$\text{PI}^{(p)}(\beta_p)$}, $u$ must be constant. Since $u\in C_c(X)$ and $X$ is unbounded, it follows that $u=0$, contradicting the assumption that $u$ is non-zero.\qed

\section{Proof of Theorem \ref{thm_gammap}\ref{item_gammap1}}\label{sec_gammap1}

\begin{lemma}\label{lem_betap_strict}
Assume \ref{eq_VD}. Let $I\subseteq(1,+\infty)$ be an open interval. Assume that, for any $p\in I$, there exists a $p$-energy $(\mathcal{E}^{(p)},\mathcal{F}^{(p)})$ satisfying \hyperref[eq_PI]{$\text{PI}^{(p)}(\beta_p)$} and \hyperref[eq_ucap]{$\text{cap}^{(p)}(\beta_p)_\le$}. Then the function $p\mapsto \beta_p$ is \emph{strictly} increasing on $I$.
\end{lemma}

\begin{proof}
By Lemmas \ref{lem_Bau} and \ref{lem_Shi1}, the function $p\mapsto\beta_p$ is continuous and increasing on $I$. It remains to show that the monotonicity is strict. Indeed, suppose that there exist $p,q\in I$ with $p<q$ such that $\beta_p=\beta_q$. By regularity, there exists a non-zero $u\in \mathcal{F}^{(p)}\cap C_c(X)$. By the uniform continuity of $u\in C_c(X)$, for any $\varepsilon>0$, there exists $\delta>0$ such that for any $r\in(0,\delta)$ and any $x,y\in X$ with $d(x,y)<r$, $\lvert u(x)-u(y)\rvert<\varepsilon$, then
\begin{align*}
&E^{(q)}_{\beta_q,X}(u,r)=\frac{1}{r^{\beta_q}}\int_X \left(\frac{1}{V(x,r)}\int_{B(x,r)}\lvert u(x)-u(y)\rvert^q m(\mathrm{d}y)\right)m(\mathrm{d}x)\\
&\le \frac{\varepsilon^{q-p}}{r^{\beta_p}}\int_X \left(\frac{1}{V(x,r)}\int_{B(x,r)}\lvert u(x)-u(y)\rvert^p m(\mathrm{d}y)\right)m(\mathrm{d}x)=\varepsilon^{q-p}E^{(p)}_{\beta_p,X}(u,r).
\end{align*}
By Lemma \ref{lem_Shi3},
$$\liminf_{r\downarrow0}E^{(q)}_{\beta_q,X}(u,r)\le \varepsilon^{q-p}\liminf_{r\downarrow0}E^{(p)}_{\beta_p,X}(u,r)\lesssim \varepsilon^{q-p}\mathcal{E}^{(p)}(u)<+\infty,$$
which gives $u\in \mathcal{F}^{(q)}$ and $\mathcal{E}^{(q)}(u)\lesssim \varepsilon^{q-p}\mathcal{E}^{(p)}(u)$. Since $\varepsilon>0$ is arbitrary, $\mathcal{E}^{(q)}(u)=0$. By \hyperref[eq_PI]{$\text{PI}^{(q)}(\beta_q)$}, $u$ must be constant. Since $u\in C_c(X)$ and $X$ is unbounded, it follows that $u=0$, contradicting the assumption that $u$ is non-zero.
\end{proof}

Fix $p,q\in I$ with $p<q$. If $\beta_p\le d_h$, then by Lemma \ref{lem_betap_strict},
\begin{equation}\label{eq_easy_strict}
\gamma_p=\frac{\beta_p-d_h}{p-1}\le \frac{\beta_p-d_h}{q-1}< \frac{\beta_q-d_h}{q-1}=\gamma_q.
\end{equation}
We may assume that $\beta_q>\beta_p>d_h$. It suffices to prove the following result.

\begin{lemma}\label{lem_u_pq}
Assume \ref{eq_Vdh} and $\beta_p>d_h$. Let $(\mathcal{E}^{(p)},\mathcal{F}^{(p)})$ be a $p$-energy satisfying \hyperref[eq_PI]{$\text{PI}^{(p)}(\beta_p)$} and \hyperref[eq_ucap]{$\text{cap}^{(p)}(\beta_p)_\le$}. For any ball $B=B(x_0,R)$, let $u\in \mathcal{F}^{(p)}\subseteq C(X)$ be the $\mathcal{E}^{(p)}$-capacitary potential for $(\overline{B},X\backslash 2B)$ as in Lemma \ref{lem_rec}\ref{item_pot}. Then for any $q>p$ and any $\beta\in(0,(q-1)\gamma_p+d_h]$, we have $u\in B^{q,\frac{\beta}{q}}(X)$.
\end{lemma}

\begin{proof}[Proof of Theorem \ref{thm_gammap}\ref{item_gammap1} using Lemma \ref{lem_u_pq}]
Assume $\beta_q>\beta_p>d_h$. Under \ref{eq_Vdh}, \hyperref[eq_PI]{$\text{PI}^{(p)}(\beta_p)$}, and \hyperref[eq_ucap]{$\text{cap}^{(p)}(\beta_p)_\le$}, let $u\in \mathcal{F}^{(p)}\subseteq C(X)$ be given as in Lemma \ref{lem_u_pq}, then $u$ is non-zero and $u\in B^{q,\frac{\beta}{q}}(X)$ for any $\beta\in(0,(q-1)\gamma_p+d_h]$. Under \ref{eq_Vdh}, \hyperref[eq_PI]{$\text{PI}^{(q)}(\beta_q)$}, and \hyperref[eq_ucap]{$\text{cap}^{(q)}(\beta_q)_\le$}, by Lemma \ref{lem_Shi1}, the existence of such a function $u$ would imply
$$\gamma_q=\frac{q\alpha_q(X)-d_h}{q-1}\ge \frac{(q-1)\gamma_p+d_h-d_h}{q-1}=\gamma_p.$$
\end{proof}

\begin{proof}[Proof of Lemma \ref{lem_u_pq}]
By Lemma \ref{lem_rec}, $u\in \mathcal{F}^{(p)}\subseteq C(X)$ satisfies $0\le u\le 1$ on $X$, $u=1$ on $\overline{B}$, $u=0$ on $X\backslash 2B$, $u$ is harmonic in $U:=2B\backslash \overline{B}$, and $\mathcal{E}^{(p)}(u)=\mathrm{cap}^{(p)}(\overline{B},X\backslash 2B)\asymp \frac{1}{R^{\beta_p-d_h}}$, where the last inequality follows from \hyperref[eq_ucap]{$\text{cap}^{(p)}(\beta_p)_\le$}. Moreover,
\begin{equation}\label{eq_uholder_R}
\lvert u(x)-u(y)\rvert\le C_1 \left(\frac{d(x,y)}{R}\right)^{\gamma_p}=C_2d(x,y)^{\gamma_p}\text{ for any }x,y\in X,
\end{equation}
where $C_1>0$ is a constant and $C_2=\frac{C_1}{R^{\gamma_p}}$.

For any $\beta>0$ and any $r>0$, by \ref{eq_Vdh}, \cite[Lemma 3.1]{Bau24} gives
$$E^{(q)}_{\beta,X}(u,r)\le \frac{2^{q-1}(2^{d_h}C_R^2+1)}{r^\beta}\lVert u\rVert^q_{L^q(X;m)}=\frac{C_3}{r^\beta}\lVert u\rVert^q_{L^q(X;m)},$$
where $C_3=2^{q-1}(2^{d_h}C_R^2+1)$, hence
\begin{equation}\label{eq_Eqlarger}
\sup_{r\in[R,+\infty)}E^{(q)}_{\beta,X}(u,r)\le \frac{C_3}{R^\beta}\lVert u\rVert^q_{L^q(X;m)}<+\infty.
\end{equation}
For any $r\in(0,R)$, by (\ref{eq_uholder_R}),
\begin{align*}
&E^{(q)}_{\beta,X}(u,r)=E^{(q)}_{\beta,4B}(u,r)\\
&\le \frac{1}{r^\beta}\int_{4B}\left(\frac{1}{V(x,r)}\int_{B(x,r)}\left(C_2d(x,y)^{\gamma_p}\right)^{q-p}\lvert u(x)-u(y)\rvert^p m(\mathrm{d}y)\right)m(\mathrm{d}x)\\
&\le C_2^{q-p} r^{(q-p)\gamma_p+\beta_p-\beta}E^{(p)}_{\beta_p,4B}(u,r)=C_2^{q-p} r^{(q-1)\gamma_p+d_h-\beta}E^{(p)}_{\beta_p,X}(u,r).
\end{align*}
Since $u\in \mathcal{F}^{(p)}$, by \ref{eq_Vdh}, \hyperref[eq_PI]{$\text{PI}^{(p)}(\beta_p)$}, and \hyperref[eq_ucap]{$\text{cap}^{(p)}(\beta_p)_\le$}, Lemma \ref{lem_Shi3} gives
$$\sup_{r\in(0,R)}E^{(p)}_{\beta_p,X}(u,r)<+\infty.$$
Hence for any $\beta\in(0,(q-1)\gamma_p+d_h]$,
$$\sup_{r\in(0,R)}E^{(q)}_{\beta,X}(u,r)<+\infty.$$
Combining this with (\ref{eq_Eqlarger}), $u\in B^{q,\frac{\beta}{q}}(X)$ for any $\beta\in(0,(q-1)\gamma_p+d_h]$.
\end{proof}

\section{Proof of Theorem \ref{thm_gammap}\ref{item_gammap2}}\label{sec_gammap2}

By \ref{eq_CC} and \cite[PROPOSITION A.1]{KM20}, $d$ is Lipschitz equivalent to a geodesic metric. Throughout this section, we may assume without loss of generality that $d$ is a geodesic metric.

By Lemmas \ref{lem_Bau} and \ref{lem_Shi1}, together with Theorem \ref{thm_gammap}\ref{item_gammap1} and (\ref{eq_betap_bdd}), the function $p\mapsto\gamma_p=\frac{\beta_p-d_h}{p-1}$ takes values in $(-\infty,1]$ and is continuous and increasing on $I$. Hence,
$$
J:=\{p\in I:\gamma_p=1\}
$$
is relatively closed in $I$. Moreover, if $p\in J$, then $[p,+\infty)\cap I\subseteq J$. It therefore suffices to show that, for any $p\in J$, there exists $\varepsilon>0$ such that $(p-\varepsilon,p]\subseteq J$. This would imply that $J$ is relatively open in $I$. Since $I$ is an open interval, either $J=\emptyset$ or $J=I$. In what follows, we assume that
$$
I_0:=\{p\in I:\beta_p>d_h\}\ne\emptyset
$$
and restrict our attention to parameters in $I_0$, where Lemma \ref{lem_rec} applies.

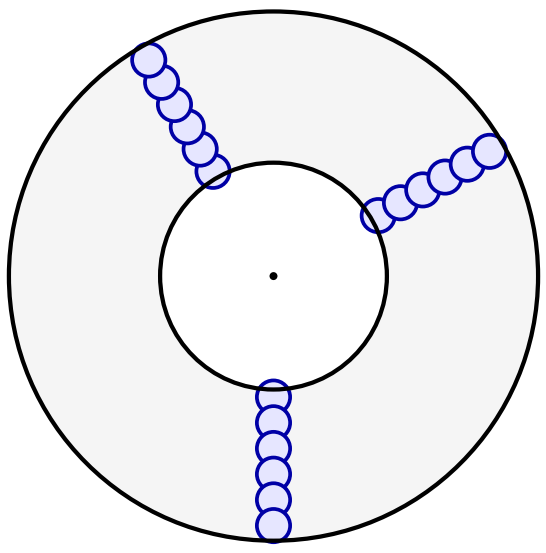
\begin{figure}[ht]
  \centering
  \begin{tikzpicture}[scale=0.5]
    \def\rin{1.5}
    \def\rout{3.5}
    \def\rsmall{0.22}

    \path[fill=gray!8,even odd rule]
      (0,0) circle (\rout)
      (0,0) circle (\rin);

    \foreach \angle in {30,120,270}{
      \foreach \radius in {1.60,1.94,2.28,2.62,2.96,3.30}{
        \draw[blue!65!black,fill=blue!10,line width=0.7pt]
          (\angle:\radius) circle (\rsmall);
      }
    }

    \draw[thick] (0,0) circle (\rin);
    \draw[thick] (0,0) circle (\rout);

    \fill (0,0) circle (1.5pt);
  \end{tikzpicture}
\caption{Disjoint chains across an annulus}\label{fig_chain}
\end{figure}

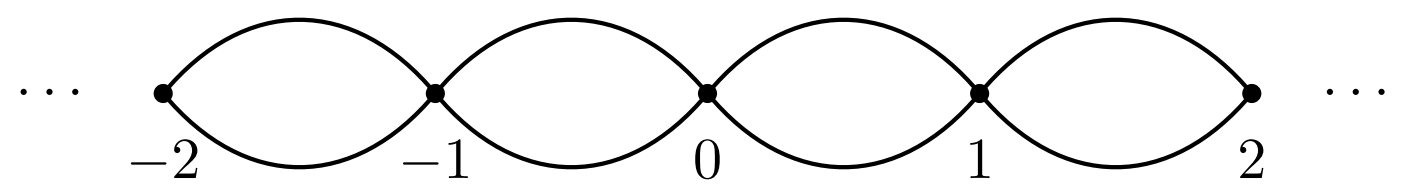
\begin{figure}[ht]
  \centering
  \begin{tikzpicture}[x=1.8cm,y=1cm]
    \foreach \i in {-2,-1,0,1}{
      \pgfmathtruncatemacro{\j}{\i+1}
      \draw[thick]
        (\i,0) .. controls (\i+0.3,0.65) and (\i+0.7,0.65) .. (\j,0);
      \draw[thick]
        (\i,0) .. controls (\i+0.3,-0.65) and (\i+0.7,-0.65) .. (\j,0);
    }

    \foreach \i in {-2,-1,0,1,2}{
      \fill (\i,0) circle (1.8pt);
      \node[below=5pt] at (\i,0) {$\i$};
    }

    \node at (-2.4,0) {$\cdots$};
    \node at (2.4,0) {$\cdots$};
  \end{tikzpicture}
  \caption{The integer lattice with doubled edges}
  \label{fig_lattice}
\end{figure}

The proof is motivated by the argument for the upper bound in \cite[Proposition 2.1]{Yan25a}, particularly by the conditions under which equality can hold. That argument suggests that, if there exists $p$ with $\gamma_p=1$, that is, $\beta_p=d_h+p-1$, then, roughly speaking, there exists an integer $N$ such that for any pair of concentric balls, at most $N$ ``disjoint" chains of small balls can connect their boundaries, as on the Vicsek set; see Figure \ref{fig_chain}. This suggests that the space should have some form of tree structure. However, connecting each pair of adjacent integers in $\mathbb{Z}$ by two parallel edges gives a space containing cycles, as illustrated in Figure \ref{fig_lattice}, but this space has the same Hausdorff dimension $d_h=1$ and the same $p$-walk dimension $\beta_p=p$ as $\mathbb{R}$. Thus, $\beta_p=d_h+p-1$. Nevertheless, this example has a local tree structure, suggesting that such a structure may still be expected at small scales. This motivates the following definition of a balanced triple, which is intended to capture this structure quantitatively.

For $r>0$ and $x,y,z\in X$, we say that $(x,y,z)$ is a balanced triple at scale $r$ if
$$\frac{1}{2}r \le \min\{d(x,y),d(x,z),d(y,z)\}\le\max\{d(x,y),d(x,z),d(y,z)\}\le2r.$$
We define its deficit by
$$\Delta^{(p)}(x,y,z;r):=\frac{1}{r^{\gamma_p}}\left(\rho^{(p)}(x,y)+\rho^{(p)}(y,z)-\rho^{(p)}(x,z)\right).$$
Then by the triangle inequality for $\rho^{(p)}$ and (\ref{eq_rhodc}), $0\le \Delta^{(p)}(x,y,z;r)\lesssim 1$. The quantity $\Delta^{(p)}$ is a normalized ordered triangle inequality deficit for the resistance metric $\rho^{(p)}$. Such deficits and triples with comparable side lengths have been considered previously; see, for example, \cite[Section 1.1, Remark 1.2]{Sch07} and \cite[Eqs. (2.2)--(2.4)]{DS21}.

The key ingredient in the proof is the following result.

\begin{proposition}\label{prop_Deltapq}
Let $p,q\in I_0$. Assume $\gamma_p=1$. Then $\Delta^{(p)}\to0$ implies $\Delta^{(q)}\to0$. More precisely, for any $\xi>0$, there exists $\kappa>0$ such that, for any $r>0$ and any balanced triple $(x,y,z)$ at scale $r$ satisfying $\Delta^{(p)}(x,y,z;r)<\kappa$, we have $\Delta^{(q)}(x,y,z;r)<\xi$.
\end{proposition}

\begin{proof}[Proof of Theorem \ref{thm_gammap}\ref{item_gammap2} using Proposition \ref{prop_Deltapq}]
For simplicity, take $x,y\in X$ with $d(x,y)=1$, then there exists a family of points
$$\left\{x_{n,k}:k=0,1,\ldots,2^n,n=0,1,\ldots\right\},$$
such that, for any $n\ge0$, $x_{n,0}=x$ and $x_{n,2^n}=y$, for any $k=0,1,\ldots,2^{n}-1$, $d(x_{n,k},x_{n,k+1})=\frac{1}{2^n}$, and $x_{n,k}=x_{n+1,2k}$. For any $s\in I_0$ and any $n\ge0$, let
$$L^{(s)}_n:=\sum_{k=0}^{2^n-1}\rho^{(s)}(x_{n,k},x_{n,k+1}).$$
Since $\rho^{(s)}(x,y)\asymp d(x,y)^{\gamma_s}$,
\begin{equation}\label{eq_dic1}
L^{(s)}_n\asymp\sum_{k=0}^{2^n-1}d(x_{n,k},x_{n,k+1})^{\gamma_s}=2^{(1-\gamma_s)n}.
\end{equation}
By the triangle inequality for $\rho^{(s)}$, $\{L^{(s)}_n\}_n$ is increasing. Moreover,
\begin{align*}
&L^{(s)}_{n+1}-L^{(s)}_n\\
&=\sum_{k=0}^{2^n-1}\left(\rho^{(s)}(x_{n+1,2k},x_{n+1,2k+1})+\rho^{(s)}(x_{n+1,2k+1},x_{n+1,2k+2})-\rho^{(s)}(x_{n,k},x_{n,k+1})\right)\\
&=\sum_{k=0}^{2^n-1}\left(\rho^{(s)}(x_{n,k},x_{n+1,2k+1})+\rho^{(s)}(x_{n+1,2k+1},x_{n,k+1})-\rho^{(s)}(x_{n,k},x_{n,k+1})\right)\\
&=\frac{1}{2^{\gamma_s n}}\sum_{k=0}^{2^n-1}\Delta^{(s)}\left(x_{n,k},x_{n+1,2k+1},x_{n,k+1};\frac{1}{2^n}\right)=\frac{1}{2^{\gamma_s n}}\sum_{k=0}^{2^n-1}\Delta^{(s)}_{n,k},
\end{align*}
where $\Delta^{(s)}_{n,k}:=\Delta^{(s)}\left(x_{n,k},x_{n+1,2k+1},x_{n,k+1};\frac{1}{2^n}\right)$ is bounded by a constant $C^{(s)}>0$ independent of $n$ and $k$. Let
$$M^{(s)}_n:=\frac{1}{2^{n}}\sum_{k=0}^{2^n-1}\Delta^{(s)}_{n,k},$$
then
\begin{equation}\label{eq_dic2}
L^{(s)}_{n+1}-L^{(s)}_n=2^{(1-\gamma_s)n}M^{(s)}_n.
\end{equation}

Since $\gamma_p=1$, by (\ref{eq_dic1}), $\{L^{(p)}_n\}_n$ is bounded and increasing. By (\ref{eq_dic2}), $\sum_{n=0}^{+\infty}M^{(p)}_n<+\infty$, which gives $\lim_{n\to+\infty}M^{(p)}_n=0$.

Let $q\in I_0$. For any $\xi>0$, by Proposition \ref{prop_Deltapq}, there exists $\kappa>0$ such that $\Delta^{(p)}_{n,k}<\kappa$ implies $\Delta^{(q)}_{n,k}<\xi$. Noting that
$$M^{(p)}_n\ge\frac{1}{2^n}\sum_{k:\Delta^{(p)}_{n,k}\ge\kappa}\Delta^{(p)}_{n,k}\ge \frac{\kappa}{2^n}\#\{k:\Delta^{(p)}_{n,k}\ge\kappa\},$$
we have
\begin{align*}
&M^{(q)}_n=\frac{1}{2^n}\left(\sum_{k:\Delta^{(p)}_{n,k}<\kappa}+\sum_{k:\Delta^{(p)}_{n,k}\ge\kappa}\right)\Delta^{(q)}_{n,k}\le \frac{1}{2^n}\sum_{k:\Delta^{(p)}_{n,k}<\kappa}\xi+\frac{1}{2^n}\sum_{k:\Delta^{(p)}_{n,k}\ge\kappa}C^{(q)}\\
&\le \xi+\frac{C^{(q)}}{2^n}\#\{k:\Delta^{(p)}_{n,k}\ge\kappa\}\le \xi+\frac{C^{(q)}}{\kappa}M^{(p)}_n.
\end{align*}
Letting $n\to+\infty$ and then $\xi\downarrow0$, we obtain $\lim_{n\to+\infty}M^{(q)}_n=0$. By (\ref{eq_dic1}) and (\ref{eq_dic2}),
$$0\le \frac{L^{(q)}_{n+1}}{L^{(q)}_n}-1=\frac{2^{(1-\gamma_q)n}}{L^{(q)}_n}M^{(q)}_n\asymp M^{(q)}_n\to0.$$
Thus $\lim_{n\to+\infty}\log \frac{L^{(q)}_{n+1}}{L^{(q)}_n}=0$, and $\lim_{N\to+\infty}\frac{1}{N}\sum_{n=0}^{N-1}\log \frac{L^{(q)}_{n+1}}{L^{(q)}_n}=0$, which gives
$$\lim_{N\to+\infty}\frac{1}{N}{\log L^{(q)}_N}=0.$$
On the other hand, by (\ref{eq_dic1}),
$$\lim_{N\to+\infty}\frac{1}{N}{\log L^{(q)}_N}=(1-\gamma_q)\log2.$$
Therefore, $\gamma_q=1$.
\end{proof}

A map $\theta:[0,1]\to X$ is called a curve if it is continuous. For $A\subseteq X$, write $\theta\cap A=\theta([0,1])\cap A$. Let $K_0,K_1,U\subseteq X$ be such that $K_0$ and $K_1$ are compact sets, $U$ is open, and $K_0\cup K_1\subseteq U$. Define
$$
\mathrm{Path}(K_0,K_1;U)
:=
\left\{\theta:[0,1]\to X:\theta\text{ is a curve},\ \theta(0)\in K_0,\ \theta(1)\in K_1,\theta([0,1])\subseteq U\right\}.
$$
For $x,y\in U$, write
$$\mathrm{Path}(x,y;U)=\mathrm{Path}(\{x\},\{y\};U).$$

The proof of Proposition \ref{prop_Deltapq} reduces to the following two results.

\begin{lemma}\label{lem_bypass1}
Let $p\in I_0$. For any $\xi>0$, there exist $\delta_0,L_0>0$ such that, for any $\delta\in(0,\delta_0]$, $L\in[L_0,+\infty)$, and $r>0$, and any balanced triple $(x,y,z)$ at scale $r$ satisfying $\Delta^{(p)}(x,y,z;r)\ge\xi$, there exists $\theta\in\mathrm{Path}(x,z;B(y,Lr))$ such that $\theta\cap B(y,\delta r)=\emptyset$.
\end{lemma}

\begin{lemma}\label{lem_bypass2}
Let $p\in I_0$. Assume $\gamma_p=1$. For any $\delta,L>0$, there exists $\kappa>0$ such that, for any $r>0$ and any balanced triple $(x,y,z)$ at scale $r$, the following holds: if there exists $\theta\in\mathrm{Path}(x,z;B(y,Lr))$ such that $\theta\cap B(y,\delta r)=\emptyset$, then $\Delta^{(p)}(x,y,z;r)\ge\kappa$.
\end{lemma}

\begin{proof}[Proof of Proposition \ref{prop_Deltapq} using Lemmas \ref{lem_bypass1} and \ref{lem_bypass2}]
Let $p,q\in I_0$. Assume $\gamma_p=1$. Let $(x,y,z)$ be a balanced triple at scale $r$. For any $\xi>0$, choose $\delta_0,L_0>0$ as in Lemma \ref{lem_bypass1} with $q$ in place of $p$, and set $\delta=\delta_0$ and $L=L_0$. For these choices, let $\kappa>0$ be as in Lemma \ref{lem_bypass2} for the exponent $p$.

We claim that if $\Delta^{(p)}(x,y,z;r)<\kappa$, then $\Delta^{(q)}(x,y,z;r)<\xi$. Indeed, suppose that $\Delta^{(q)}(x,y,z;r)\ge\xi$. By Lemma \ref{lem_bypass1}, there exists $\theta\in\mathrm{Path}(x,z;B(y,Lr))$ such that $\theta\cap B(y,\delta r)=\emptyset$. It follows from Lemma \ref{lem_bypass2} that $\Delta^{(p)}(x,y,z;r)\ge \kappa$, a contradiction.
\end{proof}

We will frequently use the following elementary identity: for any $A,B>0$,
\begin{align}
&\min\left\{
A(1-\lambda)^p+B\lambda^p:\lambda\in[0,1]
\right\}
\nonumber\\
&=
\left.
\left(A(1-\lambda)^p+B\lambda^p\right)
\vphantom{
\frac{A^{\frac{1}{p-1}}}
{A^{\frac{1}{p-1}}+B^{\frac{1}{p-1}}}
}
\right|_{\displaystyle
\lambda=
\frac{A^{\frac{1}{p-1}}}
{A^{\frac{1}{p-1}}+B^{\frac{1}{p-1}}}}
=
\frac{1}{\left(\frac{1}{A^{\frac{1}{p-1}}}+\frac{1}{B^{\frac{1}{p-1}}}\right)^{p-1}}.\label{eq_minAB}
\end{align}
The value of $\lambda$ displayed above is the unique minimizer.

For $p\in I_0$, define
\begin{equation}\label{eq_Fpe}
\mathcal{F}^{(p)}_e
:=
\left\{
u:
\begin{array}{l}
u\text{ is an }m\text{-measurable function that is finite }m\text{-a.e.,}\\
\text{and there exists an }\mathcal{E}^{(p)}\text{-Cauchy sequence }
\{u_n\}\subseteq\mathcal{F}^{(p)},\\
\text{such that }u_n\to u\ m\text{-a.e.}
\end{array}
\right\}.
\end{equation}

Arguing as in \cite[Lemma 8.6]{GY26}, we deduce from the Morrey--Sobolev inequality and an approximation argument that any $f\in \mathcal{F}_e^{(p)}$ admits a continuous version, which we still denote by $f$. The same arguments show that $(\mathcal{E}^{(p)},\mathcal{F}^{(p)}_e)$ satisfies $(\mathrm{RF1})_p$--$(\mathrm{RF4})_p$ of \cite[Definition 6.1]{KS24aV5}. The Markovian property, $p$-Clarkson's inequality, and strong sub-additivity extend to $\mathcal{F}_e^{(p)}$ by approximation. In particular, $1\in \mathcal{F}^{(p)}_e$. Moreover, for any $u\in\mathcal{F}^{(p)}_e$ and any $x,y\in X$,
$$\lvert u(x)-u(y)\rvert^p \le R^{(p)}(x,y)\mathcal{E}^{(p)}(u).$$
By approximation, closedness, and \ref{eq_Alg}, we have $\eta f\in \mathcal{F}^{(p)}\cap C_c(X)$ for any $f\in \mathcal{F}_e^{(p)}$ and any $\eta\in \mathcal{F}^{(p)}\cap C_c(X)$; see also the localization argument in the proof of \cite[Lemma 8.11]{GY26}. We can also extend the energy measures to $\mathcal{F}_e^{(p)}$, preserving their locality property.

For any non-empty disjoint compact subsets $A,B\subseteq X$, there exists a unique $u\in \mathcal{F}^{(p)}_e$ with $u=0$ on $A$ and $u=1$ on $B$, such that
$$\mathcal{E}^{(p)}(u)=\mathrm{cap}^{(p)}(A,B)=\inf \left\{\mathcal{E}^{(p)}(u):u\in \mathcal{F}^{(p)}_e,u=0\text{ on }A,u=1\text{ on }B\right\}.$$
Such a function $u$ is called the equilibrium potential between $A$ and $B$. It satisfies $0\le u\le1$ on $X$ and is harmonic in $X\backslash(A\cup B)$.

Finally, by \cite[Theorem 2.3]{Yan25d}, $(\mathcal{E}^{(p)},\mathcal{F}^{(p)})$ satisfies the elliptic Harnack inequality, which also holds for bounded non-negative harmonic functions in $\mathcal{F}_e^{(p)}$ by localization.

First, we prove Lemma \ref{lem_bypass1}, which follows directly from the following result.

\begin{lemma}\label{lem_bypass1a}
Let $p\in I_0$. There exist $\delta_0,L_0,C>0$ such that, for any $\delta\in(0,\delta_0]$, $L\in[L_0,+\infty)$, and $r>0$, and any balanced triple $(x_0,y_0,z_0)$ at scale $r$, if for any $\theta\in\mathrm{Path}(x_0,z_0;B(y_0,Lr))$, we have $\theta\cap \overline{B(y_0,2\delta r)}\ne\emptyset$, then
$$\Delta^{(p)}(x_0,y_0,z_0;r)\le C\left(\delta^{\frac{\beta_p-d_h}{p}}+\frac{1}{{L^{\beta_p-d_h}}}\right).$$
\end{lemma}

\begin{proof}[Proof of Lemma \ref{lem_bypass1} using Lemma \ref{lem_bypass1a}]
Let $\delta_0,L_0,C>0$ be as in Lemma \ref{lem_bypass1a}. For any $\xi>0$, take $\delta_1<{\delta_0}$ and $L_1>L_0$ such that
\begin{equation}\label{eq_bypass_Cxi}
C \left(\delta_1^{\frac{\beta_p-d_h}{p}}+\frac{1}{{L_1^{\beta_p-d_h}}}\right)<\xi.
\end{equation}
For any $\delta\in(0,\delta_1]$, $L\in[L_1,+\infty)$, any $r>0$, and any balanced triple $(x,y,z)$ at scale $r$ satisfying $\Delta^{(p)}(x,y,z;r)\ge\xi$, we claim that there exists $\theta\in\mathrm{Path}(x,z;B(y,Lr))$ such that $\theta\cap B(y,\delta r)=\emptyset$. Indeed, suppose that no such path exists, then any path $\theta\in\mathrm{Path}(x,z;B(y,Lr))$ satisfies $\theta\cap \overline{B(y,2\delta r)}\ne\emptyset$. By Lemma \ref{lem_bypass1a} and (\ref{eq_bypass_Cxi}), $\Delta^{(p)}(x,y,z;r)<\xi$, a contradiction.
\end{proof}

\begin{proof}[Proof of Lemma \ref{lem_bypass1a}]
For notational convenience, we omit the superscript $(p)$ when no ambiguity can arise. Let $\delta\ll1$, $L\gg1$ be chosen later, and let $U=B(y_0,Lr)$, $K=\overline{B(y_0,2\delta r)}$.

First, we show that there exists a constant $C_1>0$ such that
\begin{align}
&\rho(x_0,y_0)\ge R(x_0,\overline{B(y_0,3\delta r)})^{\frac{1}{p-1}}\ge (1-C_1\delta^{\frac{\beta_p-d_h}{p}})\rho(x_0,y_0),\label{eq_bpa1}\\
&\rho(z_0,y_0)\ge R(z_0,\overline{B(y_0,3\delta r)})^{\frac{1}{p-1}}\ge (1-C_1\delta^{\frac{\beta_p-d_h}{p}})\rho(z_0,y_0)\label{eq_bpa2}.
\end{align}
Take $h\in \mathcal{F}_e$ with $h(x_0)=1$ and $h(y_0)=0$, such that $\mathcal{E}(h)=\mathrm{cap}(x_0,y_0)\asymp r^{d_h-\beta_p}$. For any $y\in \overline{B(y_0,3\delta r)}$,
$$0\le h(y)=\lvert h(y)-h(y_0)\rvert\le \left(R(y,y_0)\mathcal{E}(h)\right)^{\frac{1}{p}}\lesssim \left((\delta r)^{\beta_p-d_h}(r^{d_h-\beta_p})\right)^{\frac{1}{p}}=\delta^{\frac{\beta_p-d_h}{p}},$$
hence
\begin{equation}\label{eq_bpa3}
M:=\sup_{y\in \overline{B(y_0,3\delta r)}}h(y)\lesssim \delta^{\frac{\beta_p-d_h}{p}}\in(0,1).
\end{equation}
Then the function $\frac{(h-M)_+}{1-M}\in \mathcal{F}_e$ takes the value 1 at $x_0$ and the value 0 on $\overline{B(y_0,3\delta r)}$, hence
\begin{align*}
&R(x_0,\overline{B(y_0,3\delta r)})^{-1}=\mathrm{cap}(x_0,\overline{B(y_0,3\delta r)})\\
&\le \mathcal{E}\left(\frac{(h-M)_+}{1-M}\right)\le \frac{1}{(1-M)^p}\mathcal{E}(h)=\frac{1}{(1-M)^p}\frac{1}{\rho(x_0,y_0)^{p-1}},
\end{align*}
which gives
$$R(x_0,\overline{B(y_0,3\delta r)})^{\frac{1}{p-1}}\ge (1-M)^{\frac{p}{p-1}}\rho(x_0,y_0)\ge \left(1- \frac{p}{p-1}M\right)\rho(x_0,y_0).$$
Combining this with (\ref{eq_bpa3}) and the fact that $R(x_0,\overline{B(y_0,3\delta r)})\le R(x_0,y_0)=\rho(x_0,y_0)^{p-1}$, we obtain (\ref{eq_bpa1}). Similarly, we have (\ref{eq_bpa2}).

Second, we give an upper bound for $\mathrm{cap}(x_0,z_0)$. Take $f,g\in \mathcal{F}_e$ with $0\le f,g\le 1$ on $X$, $f(x_0)=g(z_0)=1$ and $f=g=0$ on $\overline{B(y_0,3\delta r)}$, such that
$$\mathcal{E}(f)=\mathrm{cap}(x_0,\overline{B(y_0,3\delta r)}),\qquad \mathcal{E}(g)=\mathrm{cap}(z_0,\overline{B(y_0,3\delta r)}).$$
By \cite[Theorem 2.3]{Yan25c}, there exists $\phi\in \mathcal{F}\cap C_c(X)$ satisfying $0\le \phi\le 1$ on $X$, $\phi=1$ on $\overline{B(y_0,\frac{L}{2}r)}$, $\mathrm{supp}(\phi)\subseteq B(y_0,Lr)$, and $\mathcal{E}(\phi)\le \frac{C_2}{(Lr)^{\beta_p-d_h}}$, where $C_2>0$ is a constant. Then $f\wedge\phi,g\wedge \phi\in \mathcal{F}\cap C_c(X)$, $\mathrm{supp}(f\wedge \phi),\mathrm{supp}(g\wedge \phi)\subseteq B(y_0,Lr)\backslash B(y_0,3\delta r)$, $(f\wedge \phi)(x_0)=(g\wedge \phi)(z_0)=1$, $f\wedge \phi=g\wedge \phi=0$ on $\overline{B(y_0,3\delta r)}$. By \ref{eq_SubAdd},
\begin{equation}\label{eq_bpa4}
\begin{aligned}
&\mathcal{E}(f\wedge \phi)\le \mathcal{E}(f)+\mathcal{E}(\phi)\le  \mathrm{cap}(x_0,\overline{B(y_0,3\delta r)})+\frac{C_2}{(Lr)^{\beta_p-d_h}},\\
&\mathcal{E}(g\wedge \phi)\le \mathcal{E}(g)+\mathcal{E}(\phi)\le  \mathrm{cap}(z_0,\overline{B(y_0,3\delta r)})+\frac{C_2}{(Lr)^{\beta_p-d_h}},
\end{aligned}
\end{equation}

Let $\Omega_1$ and $\Omega_2$ be the connected components of $B(y_0,Lr)\backslash\overline{B(y_0,2\delta r)}$ containing $x_0$ and $z_0$, respectively. Since $d$ is geodesic, and any $\theta\in\mathrm{Path}(x_0,z_0;B(y_0,Lr))$ satisfies $\theta\cap \overline{B(y_0,2\delta r)}\ne\emptyset$, we have $\Omega_1\cap \Omega_2=\emptyset$. Let $F=1_{\Omega_1}(f\wedge \phi)$ and $G=1_{\Omega_2}(g\wedge\phi)$, then by the following Lemma \ref{lem_restriction}, $F,G\in \mathcal{F}\cap C_c(X)$, $\mathrm{d}\Gamma(F)=1_{\Omega_1}\mathrm{d}\Gamma(f\wedge\phi)$, and $\mathrm{d}\Gamma(G)=1_{\Omega_2}\mathrm{d}\Gamma(g\wedge\phi)$.

For any $t\in [0,1]$, the function $(1-t)F+t(1-G)\in \mathcal{F}_e$ takes the value 1 at $x_0$ and the value $0$ at $z_0$, by the strongly local property,
$$\mathrm{cap}(x_0,z_0)\le \mathcal{E}((1-t)F+t(1-G))=(1-t)^p \mathcal{E}(F)+t^{p}\mathcal{E}(G)\le (1-t)^p \mathcal{E}(f\wedge\phi)+t^{p}\mathcal{E}(g\wedge \phi).$$
Taking the infimum over $t\in[0,1]$ and using (\ref{eq_minAB}), we obtain
$$\mathrm{cap}(x_0,z_0)\le \frac{1}{\left(\frac{1}{\mathcal{E}(f\wedge \phi)^{\frac{1}{p-1}}}+\frac{1}{\mathcal{E}(g\wedge \phi)^{\frac{1}{p-1}}}\right)^{p-1}}.$$
Combining this with (\ref{eq_bpa1}), (\ref{eq_bpa2}), and (\ref{eq_bpa4}),
\begin{align*}
& \rho(x_0,z_0)=\frac{1}{\mathrm{cap}(x_0,z_0)^{\frac{1}{p-1}}}\\
&\ge \frac{1}{\left(\mathrm{cap}(x_0,\overline{B(y_0,3\delta r)})+\frac{C_2}{(Lr)^{\beta_p-d_h}}\right)^{\frac{1}{p-1}}}+\frac{1}{\left(\mathrm{cap}(z_0,\overline{B(y_0,3\delta r)})+\frac{C_2}{(Lr)^{\beta_p-d_h}}\right)^{\frac{1}{p-1}}}\\
&\ge (1-C_1\delta^{\frac{\beta_p-d_h}{p}})(\rho(x_0,y_0)+\rho(z_0,y_0))-C_3\frac{r^{\gamma_p}}{L^{\beta_p-d_h}},
\end{align*}
where, in the last inequality, we used the fact that $\rho(x_0,y_0),\rho(z_0,y_0)\asymp r^{\gamma_p}$, and $C_3>0$ is a constant. Hence
\begin{align*}
\Delta(x_0,y_0,z_0;r)=\frac{1}{r^{\gamma_p}} \left(\rho(x_0,y_0)+\rho(z_0,y_0)-\rho(x_0,z_0)\right)\lesssim \delta^{\frac{\beta_p-d_h}{p}}+\frac{1}{{L^{\beta_p-d_h}}}.
\end{align*}
\end{proof}

\begin{lemma}\label{lem_restriction}
Let $U$ be an open set and $K\subseteq U$ a compact set. Let $\Omega$ be a connected component of $U\backslash K$. Let $u\in \mathcal{F}^{(p)}\cap C_c(X)$ satisfy $\mathrm{supp}(u)\subseteq U$ and $u=0$ in an open neighborhood of $K$. Then $1_\Omega u\in \mathcal{F}^{(p)}\cap C_c(X)$ and $\mathrm{d}\Gamma^{(p)}(1_{\Omega}u)=1_{\Omega}\mathrm{d}\Gamma^{(p)}(u)$.
\end{lemma}

\begin{proof}
Since $u\in C_c(X)$, $\mathrm{supp}(u)\subseteq U$, and $u=0$ in an open neighborhood of $K$, $\mathrm{supp}(u)$ is a compact subset of $U\backslash K$. Since $\Omega$ is a connected component of $U\backslash K$, $\Omega$ is both relatively open and relatively closed in $U\backslash K$. Hence there exists a closed set $F$ such that $\Omega=F\cap (U\backslash K)$, which gives
$$\mathrm{supp}(u)\cap \Omega=\mathrm{supp}(u)\cap F\cap (U\backslash K)=\mathrm{supp}(u)\cap F,$$
where we used $\mathrm{supp}(u)\subseteq U\backslash K$ in the second equality. Hence $\mathrm{supp}(u)\cap \Omega$ is compact. Since $U\backslash K$ is open, $\Omega$ is open. Therefore, $\mathrm{supp}(u)\cap \Omega$ is a compact subset of the open set $\Omega$. By the regular property of $(\mathcal{E}^{(p)}, \mathcal{F}^{(p)})$, there exists $\phi\in \mathcal{F}^{(p)}\cap C_c(X)$ such that $0\le\phi\le1$ on $X$, $\phi=1$ on $\mathrm{supp}(u)\cap \Omega$, and $\mathrm{supp}(\phi)\subseteq \Omega$. By \ref{eq_Alg}, $1_\Omega u=\phi u\in \mathcal{F}^{(p)}\cap C_c(X)$. By the strongly local property, $\mathrm{d}\Gamma^{(p)}(1_\Omega u)=1_\Omega \mathrm{d}\Gamma^{(p)}(u)$.
\end{proof}

Second, we prove Lemma \ref{lem_bypass2} by contradiction.

\begin{proof}[Proof of Lemma \ref{lem_bypass2}]
For notational convenience, we omit the superscript $(p)$ when no ambiguity can arise. Since $\gamma_p=1$, $R(x,y)\asymp d(x,y)^{\beta_p-d_h}=d(x,y)^{p-1}$ and $\rho(x,y)\asymp d(x,y)$. Suppose that this result does not hold, then there exist $\delta,L>0$, a sequence $\{r_n\}\subseteq(0,+\infty)$, and a sequence of triples $\{(x_n,y_n,z_n)\}$ such that, for any $n$, $(x_n,y_n,z_n)$ is balanced at scale $r_n$ and there exists $\theta_n\in\mathrm{Path}(x_n,z_n;B(y_n,Lr_n))$ satisfying $\theta_n\cap B(y_n,\delta r_n)=\emptyset$, whereas $\Delta_n:=\Delta(x_n,y_n,z_n;r_n)\to0$.

For any $n$, let $u_n\in \mathcal{F}_e$ be the equilibrium potential between $x_n$ and $z_n$ with $u_n(x_n)=1$ and $u_n(z_n)=0$, then $\mathrm{cap}(x_n,z_n)=\mathcal{E}(u_n)$. Let $\lambda_n:=u_n(y_n)\in[0,1]$. If $\lambda_n\in(0,1)$, let
$$v_n:=\frac{(u_n-\lambda_n)_+}{1-\lambda_n},\qquad w_n:=\frac{(\lambda_n-u_n)_+}{\lambda_n},$$
then $v_n,w_n\in \mathcal{F}_e$ satisfy $v_n(x_n)=w_n(z_n)=1$ and $v_n(y_n)=w_n(y_n)=0$, hence
\begin{equation}\label{eq_bpb0}
\begin{aligned}
&\mathrm{cap}(x_n,y_n)\le\mathcal{E}(v_n)= \frac{1}{(1-\lambda_n)^p}\mathcal{E}((u_n-\lambda_n)_+),\\
&\mathrm{cap}(z_n,y_n)\le\mathcal{E}(w_n)= \frac{1}{\lambda_n^p}\mathcal{E}((\lambda_n-u_n)_+).
\end{aligned}
\end{equation}
By the strongly local property,
$$\mathrm{cap}(x_n,z_n)=\mathcal{E}(u_n)=\mathcal{E}((u_n-\lambda_n)_+)+\mathcal{E}((\lambda_n-u_n)_+)\ge (1-\lambda_n)^p \mathrm{cap}(x_n,y_n)+\lambda_n^p \mathrm{cap}(z_n,y_n),$$
which gives
$$\frac{1}{\rho(x_n,z_n)^{p-1}}\ge \frac{(1-\lambda_n)^p}{\rho(x_n,y_n)^{p-1}}+\frac{\lambda_n^p}{\rho(z_n,y_n)^{p-1}},$$
that is,
\begin{equation}\label{eq_bpb1}
1\ge {(1-\lambda_n)^p}\left(\frac{{\rho(x_n,z_n)}}{\rho(x_n,y_n)}\right)^{p-1}+{\lambda_n^p}\left(\frac{{\rho(x_n,z_n)}}{\rho(z_n,y_n)}\right)^{p-1}.
\end{equation}
If $\lambda_n=u_n(y_n)=0$, then $u_n(x_n)=1$ gives $\mathrm{cap}(x_n,y_n)\le \mathcal{E}(u_n)=\mathrm{cap}(x_n,z_n)$, that is, $\frac{\rho(x_n,z_n)}{\rho(x_n,y_n)}\le1$. Hence (\ref{eq_bpb1}) still holds. Similarly, if $\lambda_n=u_n(y_n)=1$, then (\ref{eq_bpb1}) also holds.

By (\ref{eq_minAB}),
\begin{align}
&{(1-\lambda_n)^p}\left(\frac{{\rho(x_n,z_n)}}{\rho(x_n,y_n)}\right)^{p-1}+{\lambda_n^p}\left(\frac{{\rho(x_n,z_n)}}{\rho(z_n,y_n)}\right)^{p-1}\nonumber\\
&\ge \min \left\{{(1-\lambda)^p}\left(\frac{{\rho(x_n,z_n)}}{\rho(x_n,y_n)}\right)^{p-1}+{\lambda^p}\left(\frac{{\rho(x_n,z_n)}}{\rho(z_n,y_n)}\right)^{p-1}:\lambda\in [0,1]\right\}\nonumber\\
&= \frac{1}{\left(\frac{\rho(x_n,y_n)}{\rho(x_n,z_n)}+\frac{\rho(z_n,y_n)}{\rho(x_n,z_n)}\right)^{p-1}}=\left(\frac{\rho(x_n,z_n)}{\rho(x_n,y_n)+\rho(z_n,y_n)}\right)^{p-1},\label{eq_bpb2}
\end{align}
where the minimum is attained uniquely at $\lambda=\frac{\frac{1}{\rho(x_n,y_n)}}{\frac{1}{\rho(x_n,y_n)}+\frac{1}{\rho(z_n,y_n)}}$. Recall that
$$0\le\Delta_n=\frac{1}{r_n}\left(\rho(x_n,y_n)+\rho(z_n,y_n)-\rho(x_n,z_n)\right)\to0.$$
Combining this with (\ref{eq_bpb1}) and (\ref{eq_bpb2}),
\begin{equation}\label{eq_bpb3}
1\ge \left(\frac{\rho(x_n,z_n)}{\rho(x_n,y_n)+\rho(z_n,y_n)}\right)^{p-1}=\left(1- \frac{r_n}{\rho(x_n,y_n)+\rho(z_n,y_n)}\Delta_n\right)^{p-1}.
\end{equation}
Since $\rho(x_n,y_n)\asymp d(x_n,y_n)\asymp r_n$ and $\rho(z_n,y_n)\asymp d(z_n,y_n)\asymp r_n$, $\frac{r_n}{\rho(x_n,y_n)+\rho(z_n,y_n)}\asymp1$. Since $\Delta_n\to0$, the RHS of (\ref{eq_bpb3}) converges to $1$. This has the following important consequences: First, by the continuity with respect to $\lambda$ of the function appearing in (\ref{eq_bpb2}),
$$\lambda_n-\frac{\frac{1}{\rho(x_n,y_n)}}{\frac{1}{\rho(x_n,y_n)}+\frac{1}{\rho(z_n,y_n)}}\to0,$$
where the sequence $\{\frac{\frac{1}{\rho(x_n,y_n)}}{\frac{1}{\rho(x_n,y_n)}+\frac{1}{\rho(z_n,y_n)}}\}\subseteq[0,1]$ is uniformly bounded away from 0 and 1, hence $\{\lambda_n\}\subseteq[0,1]$ is also uniformly bounded away from 0 and 1. Second, the RHS of (\ref{eq_bpb1}) converges to 1. Combining this with (\ref{eq_bpb0}) and the fact that both $\{1-\lambda_n\}$ and $\{\lambda_n\}$ are uniformly bounded away from 0, we obtain
$$1\le \frac{\mathcal{E}(v_n)}{\mathrm{cap}(x_n,y_n)}\to1,\quad 1\le \frac{\mathcal{E}(w_n)}{\mathrm{cap}(z_n,y_n)}\to1.$$
Since $\mathrm{cap}(x_n,y_n)\asymp d(x_n,y_n)^{1-p}\asymp r_n^{1-p}$ and $\mathrm{cap}(z_n,y_n)\asymp d(z_n,y_n)^{1-p}\asymp r_n^{1-p}$,
\begin{equation}\label{eq_bpb4}
\mathcal{E}(v_n)=\mathrm{cap}(x_n,y_n)+o(r_n^{1-p}),\qquad \mathcal{E}(w_n)=\mathrm{cap}(z_n,y_n)+o(r_n^{1-p}).
\end{equation}

Let $V_n\in \mathcal{F}_e$ be the equilibrium potential between $x_n$ and $y_n$ with $V_n(x_n)=1$ and $V_n(y_n)=0$, then $\mathrm{cap}(x_n,y_n)=\mathcal{E}(V_n)$. We claim that
\begin{equation}\label{eq_bpb5}
\mathcal{E}(v_n-V_n)=o(r_n^{1-p}).
\end{equation}
Indeed, if $p\in(1,2]$, then by \ref{eq_Cla},
\begin{align*}
&\mathcal{E}(v_n)+\mathrm{cap}(x_n,y_n)=\mathcal{E}(v_n)+\mathcal{E}(V_n)\\
&\ge 2 \left(\mathcal{E}\left(\frac{v_n+V_n}{2}\right)^{\frac{1}{p-1}}+\mathcal{E}\left(\frac{v_n-V_n}{2}\right)^{\frac{1}{p-1}}\right)^{p-1}\\
&\ge 2 \left(\mathrm{cap}(x_n,y_n)^{\frac{1}{p-1}}+\mathcal{E}\left(\frac{v_n-V_n}{2}\right)^{\frac{1}{p-1}}\right)^{p-1},
\end{align*}
where the last inequality follows from the fact that $\frac{v_n+V_n}{2}\in \mathcal{F}_e$ takes the value 1 at $x_n$ and the value 0 at $y_n$. By (\ref{eq_bpb4}), we obtain (\ref{eq_bpb5}). Similarly, if $p\in[2,+\infty)$, then \ref{eq_Cla} and (\ref{eq_bpb4}) also give (\ref{eq_bpb5}).

Finally, we consider $V_n(z_n)$ and derive a contradiction. On the one hand, by (\ref{eq_bpb5}),
\begin{align*}
&\lvert V_n(z_n)\rvert=\lvert (v_n-V_n)(z_n)\rvert=\lvert (v_n-V_n)(z_n)-(v_n-V_n)(y_n)\rvert\le \left(R(z_n,y_n)\mathcal{E}(v_n-V_n)\right)^{\frac{1}{p}}\\
&\asymp \left(d(z_n,y_n)^{p-1}\mathcal{E}(v_n-V_n)\right)^{\frac{1}{p}}\asymp \left(r_n^{p-1}\mathcal{E}(v_n-V_n)\right)^{\frac{1}{p}}=\left(r_n^{p-1}o(r_n^{1-p})\right)^{\frac{1}{p}}=o(1),
\end{align*}
hence $V_n(z_n)\to0$. On the other hand, recall that $\theta_n\in\mathrm{Path}(x_n,z_n;B(y_n,Lr_n))$ satisfies $\theta_n\cap B(y_n,\delta r_n)=\emptyset$. Let $\varepsilon\in(0,\frac{\delta}{2}\wedge\frac{1}{4})$ be chosen later, then for any $x\in \overline{B(x_n,\varepsilon r_n)}$,
\begin{align*}
&\lvert V_n(x)-1\rvert=\lvert V_n(x)-V_n(x_n)\rvert\le \left(R(x,x_n)\mathcal{E}(V_n)\right)^{\frac{1}{p}}\\
&\asymp \left(d(x,x_n)^{p-1}\mathrm{cap}(x_n,y_n)\right)^{\frac{1}{p}}\lesssim \left((\varepsilon r_n)^{p-1}r_n^{1-p}\right)^{\frac{1}{p}}=\varepsilon^{\frac{p-1}{p}}.
\end{align*}
Hence there exists a constant $\varepsilon\ll1$, such that $\lvert V_n(x)-1\rvert\le \frac{1}{2}$ for any $x\in \overline{B(x_n,\varepsilon r_n)}$, which gives $V_n\ge \frac{1}{2}$ on $\overline{B(x_n,\varepsilon r_n)}$; in particular, let $t_n:=\sup\{t\in[0,1]:\theta_n(t)\in \partial B(x_n,\varepsilon r_n)\}$ and $\widetilde{x}_n:=\theta_n(t_n)$, then $t_n<1$, $\{\widetilde{x}_n\}=( \theta_n|_{[t_n,1]})\cap \partial B(x_n,\varepsilon r_n)$, and $V_n(\widetilde{x}_n)\ge \frac{1}{2}$; see Figure \ref{fig_Harnack}. Since $V_n$ is harmonic in $X\backslash\{x_n,y_n\}$, we may apply the elliptic Harnack inequality along a suitable Harnack chain following $\theta_n|_{[t_n,1]}$; see \cite[Definition 2.10]{Che25}. Indeed, since $\theta_n|_{[t_n,1]}$ has distance at least $\varepsilon r_n$ from $x_n$ and at least $\delta r_n$ from $y_n$, there exists $a>0$, independent of $n$, sufficiently small such that the Harnack dilates of balls of radius $ar_n$ centered on $\theta_n|_{[t_n,1]}$ avoid $\{x_n,y_n\}$. Choose a maximal $(\frac{1}{2}ar_n)$-separated subset of $\theta_n|_{[t_n,1]}$. The balls of radius $ar_n$ centered at these points cover $\theta_n|_{[t_n,1]}$. Since $\theta_n|_{[t_n,1]}$ is connected, there exists a chain of overlapping balls in this cover joining $\widetilde{x}_n$ to $z_n$ without repeating balls. By \ref{eq_VD} and the fact that $\theta_n|_{[t_n,1]}\subseteq B(y_n,Lr_n)$, the number of balls in the covering, and hence the length of the chain, is bounded independently of $n$. Applying the elliptic Harnack inequality along this chain and using $V_n(\widetilde{x}_n)\ge \frac{1}{2}$, we obtain $V_n(z_n)\ge c>0$, where $c>0$ is a constant independent of $n$. This contradicts $V_n(z_n)\to0$.

\begin{figure}
\centering
\begin{tikzpicture}[
    x=1cm,y=1cm,
    line cap=round,line join=round,
    ball/.style={draw=black!65,fill=black!4,line width=0.6pt,dashed},
balldd/.style={draw=black!65,fill=black!4,line width=0.6pt},
    curve/.style={draw=black,line width=1.05pt},
    direction/.style={postaction={decorate},decoration={markings,
        mark=at position 0.65 with {\arrow{Stealth[length=2.0mm,width=1.4mm]}}}},
    every node/.style={font=\normalsize,inner sep=2pt}
]
    \coordinate (xn) at (0,0);
    \coordinate (yn) at (3,4.8);
    \coordinate (zn) at (6,0);
    \coordinate (xtilde) at (0.42,0);

    \draw[balldd] (xn) circle[radius=0.42];
    \draw[ball] (yn) circle[radius=1.00];

    \draw[curve,direction]
        (xn) -- (xtilde)
        .. controls (0.95,0) and (1.00,1.45) .. (1.80,2.45)
        .. controls (2.45,3.30) and (3.15,3.55) .. (3.85,2.90)
        .. controls (4.75,2.05) and (5.05,0.35) .. (zn);

    \fill (xn) circle[radius=1.45pt];
    \fill (yn) circle[radius=1.45pt];
    \fill (zn) circle[radius=1.45pt];
    \fill (xtilde) circle[radius=1.65pt];

    \node[above=0pt] at (xn) {$x_n$};
    \node[below right=4pt] at (xtilde) {$\widetilde{x}_n$};
    \node[below right=2pt] at (yn) {$y_n$};
    \node[below right=2pt] at (zn) {$z_n$};
    \node[above left] at (-0.15,0.43) {$\overline{B(x_n,\varepsilon r_n)}$};
    \node[above=5pt] at (3,5.8) {$B(y_n,\delta r_n)$};
    \node[above left] at (1.48,2.05) {$\theta_n$};
\end{tikzpicture}
\caption{The path $\theta_n$ avoiding $B(y_n,\delta r_n)$}\label{fig_Harnack}
\end{figure}
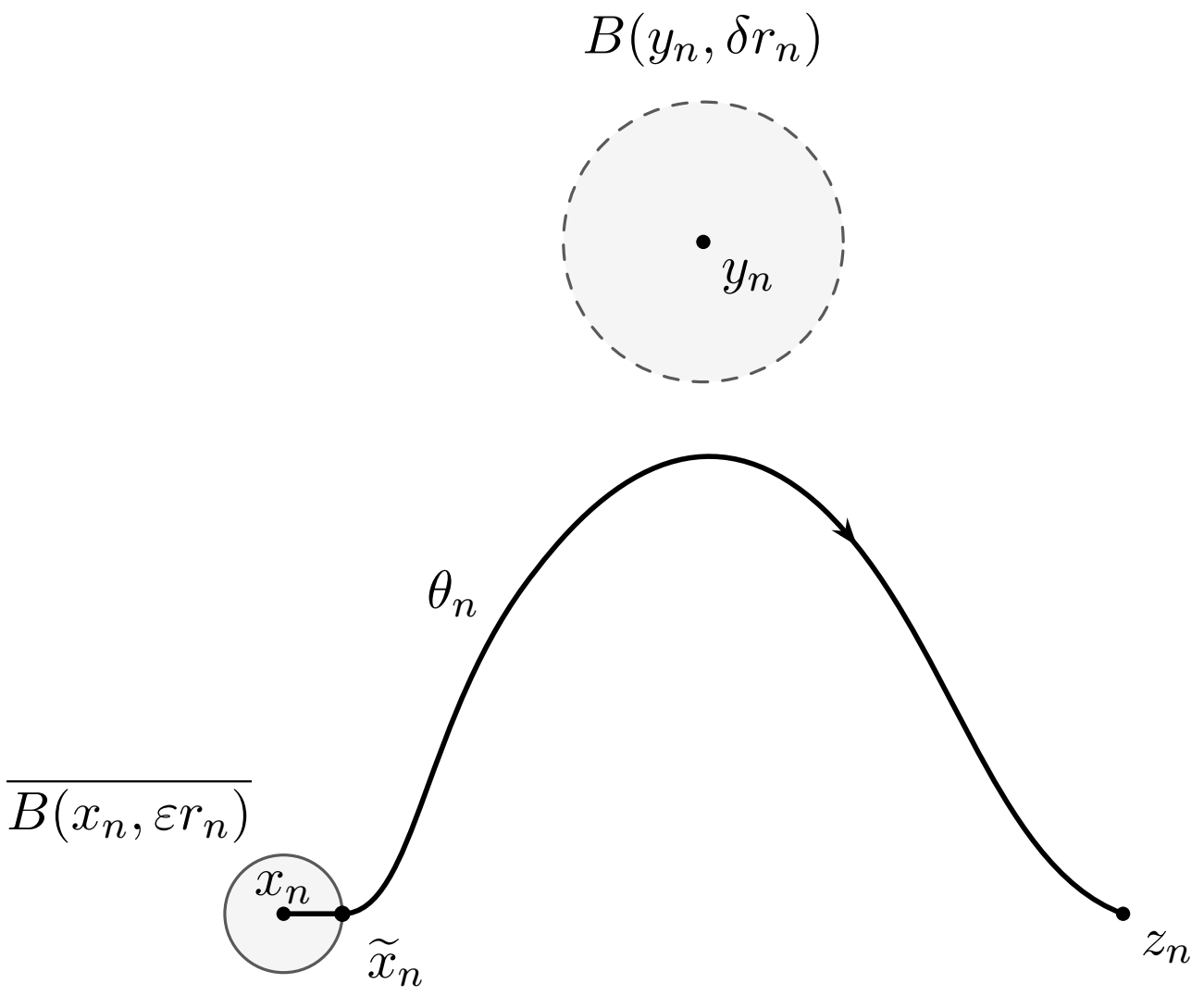
\end{proof}

\section{Proof of Theorem \ref{thm_gammap}\ref{item_gammap3}}\label{sec_gammap3}

By (\ref{eq_easy_strict}), we only need to prove the following result.

\begin{proposition}\label{prop_hard_strict}
Let $p,q\in I$ satisfy $p<q$ and $\gamma_p\in(0,1)$. Then $\gamma_p<\gamma_q$.
\end{proposition}

By Theorem \ref{thm_gammap}\ref{item_gammap1}, $\gamma_p\le\gamma_q$. Suppose that
\begin{equation}\label{eq_gammapqc}
c:=\gamma_p=\gamma_q\in(0,1).
\end{equation}
We will show that this assumption leads to a contradiction. For notational convenience, let
$$\rho:=\rho^{(p)}=(R^{(p)})^{\frac{1}{p-1}}$$
be the metric given by (\ref{eq_rhop}) associated with $(\mathcal{E}^{(p)},\mathcal{F}^{(p)})$, and write $B_\rho(x,r)$ and $\mathcal{H}^1_\rho$ for the ball and the $1$-dimensional Hausdorff measure on $(X,\rho)$.

The intuition behind the proof is as follows. First, we construct a suitable function $u$, which is in fact an $\mathcal{E}^{(p)}$-capacitary potential, such that its energy measure $\Gamma^{(p)}(u)$ admits ``good" upper and lower estimates; see Definition \ref{def_umu} and Lemma \ref{lem_Wxrn}. In particular, the equality (\ref{eq_gammapqc}) plays a key role in obtaining these lower estimates. Second, using the cutoff functions and the partition of unity from Lemma \ref{lem_cutpart}, we employ the measure $\Gamma^{(p)}(u)$ to construct a finite graph. We then apply the discrete potential theory from Proposition \ref{prop_flow}, in particular, some estimates for flows, potentials, and transportation costs, to obtain the existence of suitable paths in this graph. This enables us to show that, for any sufficiently small $\tau\ll 1$, there exist $r\ll 1$ and $\{x_{i_k}:k=0,\ldots,m\}\subseteq X$ such that
$$\rho(x_{i_0},x_{i_m})\gtrsim r,\qquad\sum_{k=1}^m\rho(x_{i_{k-1}},x_{i_k})\lesssim r,$$
and
$$\rho(x_{i_{k-1}},x_{i_k})\lesssim \tau r\text{ for any }k=1,\ldots,m.$$
Then applying (\ref{eq_rhodc}) and the triangle inequality for the metric $d$,
\begin{align*}
&r\lesssim\rho(x_{i_0},x_{i_m})\asymp d(x_{i_0},x_{i_m})^c \le \left(\sum_{k=1}^m d(x_{i_{k-1}},x_{i_k})\right)^c\asymp \left(\sum_{k=1}^m \rho(x_{i_{k-1}},x_{i_k})^{\frac{1}{c}}\right)^c\\
&\le \left(\max_{k=1,\ldots,m}\rho(x_{i_{k-1}},x_{i_k})^{\frac{1}{c}-1}\sum_{k=1}^m \rho(x_{i_{k-1}},x_{i_k})\right)^c\lesssim \left((\tau r)^{\frac{1}{c}-1}r\right)^c=\tau^{1-c}r,
\end{align*}
hence $1\lesssim \tau^{1-c}$. Since $c\in(0,1)$, letting $\tau\downarrow0$, we obtain a contradiction.

\begin{definition}\label{def_umu}
Fix a ball $B=B(x_0,R)$, let $u\in\mathcal{F}^{(p)}\subseteq C(X)$ be the $\mathcal{E}^{(p)}$-capacitary potential for $(\overline{B},X\backslash 2B)$, and define
$$\Omega:={2B}\backslash \overline{B},\qquad\mu:=\Gamma^{(p)}(u).$$
\end{definition}

By Lemma \ref{lem_rec}\ref{item_pot}, $0\le u\le 1$ on $X$, $u=1$ on $\overline{B}$, $u=0$ on $X\backslash 2B$, $u$ is harmonic in $\Omega$, $\mathcal{E}^{(p)}(u)=\mathrm{cap}^{(p)}(\overline{B},X\backslash 2B)$, and $\mu(X\backslash \Omega)=0$. Moreover, by (\ref{eq_gammapqc}), Lemmas \ref{lem_u_pq} and \ref{lem_Shi3} give that $u\in B^{q,\frac{\beta_q}{q}}(X)=\mathcal{F}^{(q)}$ and $\mathcal{E}^{(q)}(u)\in(0,+\infty)$.

\begin{lemma}\label{lem_Wxrn}
Let $u$ and $\mu$ be as in Definition \ref{def_umu}. Then there exist constants $\eta,\epsilon>0$, a Borel subset $W\subseteq \Omega$, $x\in W$, and $\{r_n\}\subseteq(0,+\infty)$ with $r_n\downarrow0$, such that $\mu(W)>0$, $\mu|_W\ge \eta\mathcal{H}^1_\rho|_W$, that is, for any Borel subset $A\subseteq W$,
\begin{equation}\label{eq_mu_lower}
\mu(A)\ge\eta \mathcal{H}^1_\rho(A),
\end{equation}
and, moreover,
\begin{equation}\label{eq_muoutW}
\lim_{r\downarrow0}\frac{\mu(B_\rho(x,r)\backslash W)}{r}=0,
\end{equation}
and
\begin{equation}\label{eq_muonW}
\mu \left(B_\rho(x,r_n)\cap W\right)\ge {\epsilon}r_n\text{ for any }n.
\end{equation}
\end{lemma}

Roughly speaking, $W$ consists of points $x$ where the function $u$ has ``large" oscillations in small balls around $x$. A suitable Morrey--Sobolev inequality then yields lower estimates for the energy measure $\mu=\Gamma^{(p)}(u)$. The equality (\ref{eq_gammapqc}) will play a key role in ensuring that $W$ has positive $\mu$-measure.

\begin{proof}
By (\ref{eq_gammapqc}), for any open subset $U\subseteq X$ and any $\theta,r>0$,
\begin{align}
&E^{(q)}_{\beta_q,U}(u,r)=\frac{1}{r^{\beta_p}}\int_U \left(\frac{1}{V(x,r)}\int_{B(x,r)}\left(\frac{\lvert u(x)-u(y)\rvert}{r^c}\right)^{q-p} \lvert u(x)-u(y)\rvert^p m(\mathrm{d}y)\right)m(\mathrm{d}x)\nonumber\\
&=\frac{1}{r^{\beta_p}} \left(\int_{\{x\in U:\mathrm{osc}_{B(x,r)}u\le\theta r^c\}}+\int_{\{x\in U:\mathrm{osc}_{B(x,r)}u>\theta r^c\}}\right) \ldots m(\mathrm{d}x)=:I_1+I_2,\label{eq_Wxrn1}
\end{align}
where
\begin{equation}\label{eq_Wxrn2}
I_1\le \theta^{q-p}E^{(p)}_{\beta_p,U}(u,r)\le \theta^{q-p}\sup_{r>0}E^{(p)}_{\beta_p,X}(u,r)=\theta^{q-p}M,
\end{equation}
where $M:=\sup_{r>0}E^{(p)}_{\beta_p,X}(u,r)\in(0,+\infty)$, which follows from $u\in \mathcal{F}^{(p)}$, $\mathcal{E}^{(p)}(u)\in(0,+\infty)$, and Lemma \ref{lem_Shi3}. By (\ref{eq_ulip}) and (\ref{eq_gammapqc}), there exists $C>0$ such that for any $x\in X$ and any $r>0$,
$$\mathrm{osc}_{B(x,r)}u\le C \mathcal{E}^{(p)}(u)^{\frac{1}{p-1}}r^c=Kr^c,$$
where $K:=C \mathcal{E}^{(p)}(u)^{\frac{1}{p-1}}\in(0,+\infty)$, hence
\begin{equation}\label{eq_Wxrn3}
I_2\le K^{q-p}E^{(p)}_{\beta_p,\{x\in U:\mathrm{osc}_{B(x,r)}u>\theta r^c\}}(u,r)=K^{q-p}E^{(p)}_{\beta_p,W_r}(u,r),
\end{equation}
where
$$W_r:= \left\{x\in U:\mathrm{osc}_{B(x,r)}u>\theta r^c\right\}.$$
Since $u\in C(X)$, the set $\{x\in X:\mathrm{osc}_{B(x,r)}u>\theta r^c\}$ is an open subset of $X$, and so is $W_r$. By Lemma \ref{lem_Shi2}, there exists $C_1>0$ such that
\begin{equation}\label{eq_Wxrn4}
E^{(p)}_{\beta_p,W_r}(u,r)\le C_1 \Gamma^{(p)}(u)\left(W_r(2A^{(p)}_{PI}r)\right)=C_1 \mu\left(W_r(2A^{(p)}_{PI}r)\right).
\end{equation}
Combining (\ref{eq_Wxrn1})--(\ref{eq_Wxrn4}),
\begin{equation}\label{eq_Wxrn5}
E^{(q)}_{\beta_q,U}(u,r)\le \theta^{q-p}M+C_1K^{q-p}\mu\left(W_r(2A^{(p)}_{PI}r)\right).
\end{equation}

Since $u\in \mathcal{F}^{(q)}$ and $\mathcal{E}^{(q)}(u)\in (0,+\infty)$, there exists an open subset $U\subseteq X$ such that $\Gamma^{(q)}(u)(U)\in(0,+\infty)$. By Lemma \ref{lem_Shi4}, there exist $C_2>0$ and $\{r_n\}\subseteq (0,+\infty)$ with $r_n\downarrow0$ such that
$$\Gamma^{(q)}(u)(U)\le C_2\liminf_{r\downarrow0}E^{(q)}_{\beta_q,U}(u,r)=C_2\lim_{n\to+\infty}E^{(q)}_{\beta_q,U}(u,r_n).$$
Combining this with (\ref{eq_Wxrn5}) and taking $\theta>0$ small enough that $\theta^{q-p}M\le \frac{1}{2C_2}\Gamma^{(q)}(u)(U)$,
$$\liminf_{n\to+\infty}\mu \left(W_{r_n}(2A^{(p)}_{PI}r_n)\right)\ge \frac{\Gamma^{(q)}(u)(U)}{2C_1C_2K^{q-p}}>0.$$
Let
$$W:=\limsup_{n\to+\infty}W_{r_n}(2A^{(p)}_{PI}r_n)=\bigcap_{N\ge1}\bigcup_{n\ge N}W_{r_n}(2A^{(p)}_{PI}r_n),$$
then $W$ is a Borel subset of $X$, and since $\mu$ is a finite measure on $X$, by Fatou's lemma,
$$\mu(W)\ge \limsup_{n\to+\infty}\mu \left(W_{r_n}(2A^{(p)}_{PI}r_n)\right)\ge\liminf_{n\to+\infty}\mu \left(W_{r_n}(2A^{(p)}_{PI}r_n)\right)\ge\frac{\Gamma^{(q)}(u)(U)}{2C_1C_2K^{q-p}}>0.$$
Since $\mu(X\backslash\Omega)=0$, replacing $W$ by $W\cap \Omega$, we may assume that $W$ is a Borel subset of $\Omega$. For any $x\in W$, by definition, for any $N\ge1$, there exist $n\ge N$ and $y\in W_{r_n}$ such that $d(x,y)<2A^{(p)}_{PI}r_n$, which gives $B(x,4A^{(p)}_{PI}r_n)\supseteq B(y,r_n)$ and
$$\mathrm{osc}_{B(x,4A^{(p)}_{PI}r_n)}u\ge \mathrm{osc}_{B(y,r_n)}u> \theta r_n^c,$$
hence
\begin{equation}\label{eq_Wxrn6}
\limsup_{r\downarrow0}\frac{\mathrm{osc}_{B(x,r)}u}{r^c}\ge \frac{\theta}{(4A^{(p)}_{PI})^c}>0.
\end{equation}
By the Morrey--Sobolev inequality implied by \hyperref[eq_PI]{$\text{PI}^{(p)}(\beta_p)$} (see, for example, \cite[Lemma 5.2]{Yan25c}) and (\ref{eq_gammapqc}), there exists $C_3>0$ such that
$$\left(\mathrm{osc}_{B(x,r)}u\right)^p\le C_3 r^{\beta_p-d_h}\Gamma^{(p)}(u)(B(x,8A^{(p)}_{PI}r))=C_3r^{(p-1)c}\mu(B(x,8A^{(p)}_{PI}r)).$$
Combining this with (\ref{eq_Wxrn6}),
$$\limsup_{r\downarrow0}\frac{\mu(B(x,r))}{r^c}\ge \frac{\theta^p}{C_3(4A^{(p)}_{PI})^{pc}(8A^{(p)}_{PI})^c}>0.$$
By (\ref{eq_rhodc}), there exists $\epsilon>0$, depending only on the constant in (\ref{eq_rhodc}) and the constant on the RHS of the preceding inequality, such that
\begin{equation}\label{eq_Wxrn7}
\limsup_{r\downarrow0}\frac{\mu(B_\rho(x,r))}{r}\ge 4\epsilon>0.
\end{equation}

For any Borel subset $A\subseteq W$, recall that $\mathcal{H}^{1}_\rho(A)=\lim_{a\downarrow0}\mathcal{H}^1_{\rho,a}(A)$ and
$$\mathcal{H}^{1}_{\rho,a}(A):=\inf \left\{\sum_{i\ge1}\mathrm{diam}_\rho(E_i):A\subseteq \bigcup_{i\ge1}E_i,\ \mathrm{diam}_\rho(E_i)\le a\right\},$$
where $\mathrm{diam}_\rho(E):=\sup\{\rho(x,y):x,y\in E\}$ for $E\subseteq X$. For any $a>0$ and any open set $G\supseteq A$, for any $x\in A$, by (\ref{eq_Wxrn7}), there exists $r_x\in(0,a)$ such that $B_\rho(x,r_x)\subseteq G$ and $\mu(B_\rho(x,r_x))\ge 2{\epsilon}r_x$. By the $5B$-covering lemma (see \cite[Theorem 1.2]{Hei01}), there exists a family of pairwise disjoint balls $\{B_\rho(x_i,r_i)\}_i$ such that $A\subseteq \cup_i B_\rho(x_i,5r_i)$, then
$$\mathcal{H}^{1}_{\rho,10a}(A)\le \sum_{i}(10r_i)\le \frac{5}{\epsilon}\sum_{i}\mu(B_\rho(x_i,r_i))\le \frac{5}{\epsilon}\mu(G).$$
Letting $a\downarrow0$, $\mathcal{H}^1_\rho(A)\le \frac{5}{\epsilon}\mu(G)$. By the outer regular property of $\mu$, $\mathcal{H}^1_\rho(A)\le \frac{5}{\epsilon}\mu(A)$. Therefore, we obtain (\ref{eq_mu_lower}) with $\eta=\frac{\epsilon}{5}$.

Since $(X,d)$ is a doubling metric space, by Assouad's embedding theorem (see \cite{Ass83,Hei01,NN12}) and (\ref{eq_rhodc}), there exist a positive integer $N$ and a constant $L\ge1$ such that $(X,\rho)$ admits an $L$-bi-Lipschitz embedding into $(\mathbb{R}^N,\lvert\cdot\rvert)$, where $\lvert \cdot\rvert$ is the standard Euclidean metric. Equivalently, there exists $\iota:X\to\mathbb{R}^N$ such that
$$\frac{1}{L}\rho(x,y)\le \lvert\iota(x)-\iota(y)\rvert \le L\rho(x,y)\text{ for any }x,y\in X.$$
For any $x\in X$ and $r>0$, let
$$U(x,r):=\iota^{-1}\left(B_{\mathbb{R}^N}(\iota(x),r)\right),$$
then
\begin{equation}\label{eq_BrhoLL2}
B_\rho(x,r)\subseteq U(x,Lr)\subseteq B_\rho(x,L^2r).
\end{equation}
By the density theorem applied to the pushforward measure $\iota_{*}\mu$ on $\mathbb{R}^N$, for $\mu$-a.e. $x\in W$,
$$\lim_{r\downarrow0}\frac{(\iota_*\mu) \left(B_{\mathbb{R}^N}(\iota(x),r)\cap \iota(W)\right)}{(\iota_*\mu) \left(B_{\mathbb{R}^N}(\iota(x),r)\right)}=1,$$
that is,
$$\lim_{r\downarrow0}\frac{\mu(U(x,r)\cap W)}{\mu(U(x,r))}=1.$$
Fix such an $x\in W$, and define $f(r):=\frac{\mu(U(x,r)\backslash W)}{\mu(U(x,r))}$, then $\lim_{r\downarrow0}f(r)=0$. By (\ref{eq_BrhoLL2}) and (\ref{eq_mu_upper}),
$$\mu \left(B_\rho(x,r)\backslash W\right)\le \mu \left(U(x,Lr)\backslash W\right)=f(Lr)\mu(U(x,Lr))\le f(Lr)\mu(B_\rho(x,L^2r))\lesssim f(Lr)L^2r.$$
On the one hand, this gives $\frac{\mu \left(B_\rho(x,r)\backslash W\right)}{r}\lesssim f(Lr)\to0$ as $r\downarrow0$, hence we obtain (\ref{eq_muoutW}). On the other hand, by (\ref{eq_Wxrn7}), there exists $\{r_n\}\subseteq(0,+\infty)$ with $r_n\downarrow0$ such that $\mu(B_\rho(x,r_n))\ge {2\epsilon}r_n$. The same estimate yields that for any sufficiently large $n$, $\mu(B_\rho(x,r_n)\backslash W)\le \epsilon r_n$, hence we obtain (\ref{eq_muonW}).
\end{proof}

We have the following cutoff functions and partition of unity on $(X,\rho)$.

\begin{lemma}\label{lem_cutpart}
There exist constants $C>0$ and $A>1$ such that the following hold:
\begin{enumerate}[label=(\arabic*),ref=(\arabic*)]
\item\label{item_cutoff} For any ball $B_\rho(x,r)$ in $(X,\rho)$, there exists $\phi\in\mathcal{F}^{(p)}$ satisfying $0\le\phi\le1$ on $X$, $\phi=1$ in $B_\rho(x,r)$, $\mathrm{supp}(\phi)\subseteq B_\rho(x,Ar)$, and
\begin{equation}\label{eq_cutoff}
\mathcal{E}^{(p)}(\phi)\le\frac{C}{r^{p-1}},\qquad \mathrm{Lip}_\rho(\phi)\le\frac{C}{r}, \qquad \Gamma^{(p)}(\phi)\le\frac{C}{r^p}\mathcal{H}^1_\rho,
\end{equation}
where $\mathrm{Lip}_\rho(\phi)$ denotes the Lipschitz constant of $\phi$ with respect to $\rho$.
\item\label{item_partition} For any $r>0$ and any $r$-net $V$ of $(X,\rho)$, there exists a family of functions $\{\phi_v\in\mathcal{F}^{(p)}:v\in V\}$ satisfying the following conditions:
\begin{itemize}
\item For any $v\in V$, $0\le\phi_v\le1$ on $X$ and $\mathrm{supp}(\phi_v)\subseteq B_\rho(v,2r)$.
\item We have $\sum_{v\in V}\phi_v=1$ on $X$.
\item For any $v\in V$,
\begin{equation}\label{eq_partition}
\mathcal{E}^{(p)}(\phi_v)\le\frac{C}{r^{p-1}},\qquad \mathrm{Lip}_\rho(\phi_v)\le\frac{C}{r},\qquad \Gamma^{(p)}(\phi_v)\le\frac{C}{r^p}\mathcal{H}^1_\rho.
\end{equation}
\end{itemize}
\end{enumerate}
\end{lemma}

\begin{proof}
Let $L\ge1$ be the constant in (\ref{eq_rhodc}). For any ball $B_\rho(x,r)$ in $(X,\rho)$,
$$B_\rho(x,r)\subseteq B(x,(Lr)^{\frac{1}{c}})\subseteq B(x,(2Lr)^{\frac{1}{c}})\subseteq B_\rho(x,2L^2r).$$
Let $\phi\in \mathcal{F}^{(p)}$ be the $\mathcal{E}^{(p)}$-capacitary potential for $(\overline{B(x,(Lr)^{\frac{1}{c}})},X\backslash B(x,(2Lr)^{\frac{1}{c}}))$, then $0\le\phi\le 1$ on $X$, $\phi=1$ in $B_\rho(x,r)$, $\mathrm{supp}(\phi)\subseteq B_\rho(x,3L^2r)$, and
$$\mathcal{E}^{(p)}(\phi)\lesssim \frac{1}{\left((Lr)^{\frac{1}{c}}\right)^{\beta_p-d_h}}\le\frac{1}{r^{p-1}}.$$
By (\ref{eq_ulip}) and (\ref{eq_mu_upperA}), $\mathrm{Lip}_\rho(\phi)\lesssim \frac{1}{r}$ and $\Gamma^{(p)}(\phi)\lesssim \frac{1}{r^p}\mathcal{H}^1_\rho$, hence we obtain (\ref{eq_cutoff}). Then for any $r>0$ and any $r$-net $V$, the family $\{\phi_v:v\in V\}$ can be constructed by the same argument as in the proof of \cite[Lemma 3.2]{Yan25a}.
\end{proof}

\begin{remark}\label{rmk_cutoff}
The same argument as in the proof of \cite[Lemma 3.5]{Yan25a} shows that, for any $A>1$, there exists a constant $C>0$ such that, for any $x\in X$ and any $r>0$, there exists a function $\phi\in\mathcal{F}^{(p)}$ satisfying the conditions in \ref{item_cutoff}.
\end{remark}

We need some results from discrete potential theory. Let $G=(V,E)$ be a finite directed graph, where $V$ is the vertex set and $E$ is the directed edge set. We write $e=(u,v)$, or equivalently $u\to v$, for a directed edge from $u$ to $v$, where $u,v\in V$ are vertices.

A path $P$ in $(V,E)$ is a finite sequence of vertices $P=(v_0,v_1,\ldots,v_n)$ such that $(v_{i-1},v_i)\in E$ for any $i=1,\ldots,n$. The path $P$ is called simple if the vertices $v_0,\ldots,v_n$ are pairwise distinct. Given $A,B\subseteq V$, we say that $P$ starts in $A$ if $v_0\in A$ and ends in $B$ if $v_n\in B$. A cycle $C$ is a finite sequence $C=(v_0,v_1,\ldots,v_n)$ such that $n\ge1$, $v_n=v_0$, $(v_{i-1},v_i)\in E$ for any $i=1,\ldots,n$, and the vertices $v_0,\ldots,v_{n-1}$ are pairwise distinct. Given a path $P=(v_0,\ldots,v_n)$, we say that an edge $e\in E$ belongs to $P$, and write $e\in P$, if $e=(v_{i-1},v_i)$ for some $i=1,\ldots,n$. For a cycle $C$, the notation $e\in C$ is defined analogously.

Let $f:E\to[0,+\infty)$ be a flow on $G$. Let $\partial f:V\to \mathbb{R}$ be given by
$$(\partial f)(u):=\sum_{v\in V:v\to u}f((v,u))-\sum_{v\in V:u\to v}f((u,v)),\ u\in V.$$
Define
$$S:=\left\{u\in V:(\partial f)(u)<0\right\},\qquad T:=\left\{u\in V:(\partial f)(u)>0\right\}.$$
Thus, $S$ and $T$ are the source set and the target set of the flow $f$, respectively.

Let $l:E\to(0,+\infty)$ be a length function on $E$. Define the transportation cost $M_l(f)$ of the flow $f$ with respect to $l$ by
$$M_l(f):=\sum_{e\in E}f(e)l(e).$$
For a path $P=(v_0,\ldots,v_n)$, define its length by
$$l(P):=\sum_{i=1}^n l((v_{i-1},v_i))=\sum_{e\in P}l(e),$$
where the sum $\sum_{e\in P}$ counts occurrences with multiplicity. The length $l(C)$ of a cycle $C$ is defined analogously.

Let $H:V\to\mathbb{R}$ be a potential on $V$. Define $dH:E\to\mathbb{R}$ by
$$(dH)(e):=H(v)-H(u),\ e=(u,v)\in E.$$
For a path $P=(v_0,\ldots,v_n)$, define the potential increment of $H$ along $P$ by
$$\Delta_P H:=H(v_n)-H(v_0)=\sum_{i=1}^n(dH)((v_{i-1},v_i))=\sum_{e\in P}(dH)(e),$$
where the sum $\sum_{e\in P}$ also counts occurrences with multiplicity. The discrete Green formula then reads
\begin{equation*}
(f,dH)_E:=\sum_{e\in E}f(e)(dH)(e)=\sum_{u\in V}(\partial f)(u)H(u)=:(\partial f,H)_V.
\end{equation*}
This identity follows from a direct calculation by changing the order of summation. See also, for example, \cite[Section 2.4, Exercise 2.11]{LP16}.

Any flow can be decomposed into a linear combination of indicator functions of simple paths and cycles as follows.

\begin{lemma}[{\cite[Proposition 1.1]{Ber91}}]\label{lem_Ber91}
Let $f:E\to[0,+\infty)$ be a flow on $G$. Then there exist a finite family $\mathcal{P}$ of simple paths starting in $S$ and ending in $T$, a finite family $\mathcal{C}$ of cycles, and positive coefficients $q_P$, $P\in\mathcal{P}$, and $r_C$, $C\in\mathcal{C}$, such that
\begin{equation*}
f=\sum_{P\in\mathcal{P}}q_P1_P+\sum_{C\in\mathcal{C}}r_C1_C.
\end{equation*}
Moreover,
\begin{equation}\label{eq_Ber2}
\sum_{P\in\mathcal{P}}q_P=-\sum_{u\in S}(\partial f)(u)=\sum_{u\in T}(\partial f)(u)=\frac{1}{2}\lVert\partial f\rVert_{\ell^1(V)},
\end{equation}
where
$$\lVert\partial f\rVert_{\ell^1(V)}:=\sum_{u\in V}\lvert(\partial f)(u)\rvert$$
is the total variation of $\partial f$, equivalently, the total boundary mass of $f$. For any length function $l:E\to(0,+\infty)$,
\begin{equation}\label{eq_Ber3}
M_l(f)=\sum_{P\in\mathcal{P}}q_P l(P)+\sum_{C\in\mathcal{C}}r_C l(C).
\end{equation}
For any potential $H:V\to \mathbb{R}$,
\begin{equation}\label{eq_Ber4}
(f,dH)_E=(\partial f,H)_V=\sum_{P\in \mathcal{P}}q_P \Delta_P H.
\end{equation}
\end{lemma}

The following result will play a key role in our argument.

\begin{proposition}\label{prop_flow}
Let $f:E\to[0,+\infty)$ be a flow on $G$, let $l:E\to(0,+\infty)$ be a length function on $E$, and let $H:V\to \mathbb{R}$ be a potential on $V$. Assume that there exist constants $C_1,C_2,C_3,C_4>0$ such that
$$\lVert\partial f\rVert_{\ell^1(V)}\le C_1,\qquad M_l(f)\le C_2,\qquad (\partial f,H)_V\ge C_3,$$
and
$$\lvert H(s)-H(t)\rvert\le C_4\text{ for any }s\in S\text{ and }t\in T.$$
Then there exists a simple path $P_*$ starting in $S$ and ending in $T$ such that
$$l(P_*)\le \frac{2C_2C_4}{C_3},\qquad \Delta_{P_*}H\ge\frac{C_3}{C_1}.$$
\end{proposition}

Roughly speaking, consider a flow $f$ through a network equipped with a length function $l$ and a potential $H$. Assume that the total boundary mass $\lVert\partial f\rVert_{\ell^1(V)}$ and the transportation cost $M_l(f)$ are bounded, that the pairing $(\partial f,H)_V$ has a positive lower bound, meaning intuitively that $H$ detects a definite amount of net transport from the source set to the target set, and that the oscillation of $H$ between the source and target sets is controlled. Then one can find a source-to-target path which is \emph{simultaneously} ``short" and along which $H$ has a definite positive ``increment".

\begin{proof}
For any $\lambda>0$, by assumption, (\ref{eq_Ber4}), and (\ref{eq_Ber2}),
\begin{align*}
&C_3\le (\partial f,H)_V=\left(\sum_{P\in \mathcal{P}:\Delta_PH< \lambda}+\sum_{P\in \mathcal{P}:\Delta_PH\ge \lambda}\right)q_P \Delta_PH\\
&\le \lambda\sum_{P\in \mathcal{P}:\Delta_PH<\lambda}q_P+C_4\sum_{P\in \mathcal{P}:\Delta_PH\ge \lambda}q_P\le \frac{1}{2}\lambda C_1+C_4\sum_{P\in \mathcal{P}:\Delta_PH\ge \lambda}q_P.
\end{align*}
Taking $\lambda=\frac{C_3}{C_1}$, we obtain
$$\sum_{P\in \mathcal{P}:\Delta_PH\ge \frac{C_3}{C_1}}q_P\ge \frac{C_3}{2C_4}>0,$$
which implies $\{P\in \mathcal{P}:\Delta_PH\ge \frac{C_3}{C_1}\}\ne\emptyset$. By assumption and (\ref{eq_Ber3}),
$$C_2\ge M_l(f)\ge\sum_{P\in \mathcal{P}:\Delta_PH\ge \frac{C_3}{C_1}}q_Pl(P),$$
hence
$$\frac{\sum_{P\in \mathcal{P}:\Delta_PH\ge \frac{C_3}{C_1}}q_Pl(P)}{\sum_{P\in \mathcal{P}:\Delta_PH\ge \frac{C_3}{C_1}}q_P}\le \frac{C_2}{\frac{C_3}{2C_4}}=\frac{2C_2C_4}{C_3},$$
which implies there exists $P_*\in \mathcal{P}$ with $\Delta_{P_*}H\ge \frac{C_3}{C_1}$ such that $l(P_*)\le \frac{2C_2C_4}{C_3}$.
\end{proof}

We now prove Proposition \ref{prop_hard_strict}, which implies Theorem \ref{thm_gammap}\ref{item_gammap3}.

\begin{proof}[Proof of Proposition \ref{prop_hard_strict}]
Let $u$ and $\mu$ be as in Definition \ref{def_umu}, and let $W,x,\{r_n\}$ be as in Lemma \ref{lem_Wxrn}. Since $x\in W\subseteq \Omega=2B\backslash \overline{B}$, we have $\overline{B_\rho(x,2r_n)}\subseteq \Omega$ for any sufficiently large $n$. Let $r=r_n$. By Lemma \ref{lem_cutpart}\ref{item_cutoff} and Remark \ref{rmk_cutoff}, there exists $\chi\in \mathcal{F}^{(p)}$ satisfying $0\le\chi\le1$ on $X$, $\chi=1$ in $B_\rho(x,r)$, $\mathrm{supp}(\chi)\subseteq B_\rho(x,2r)$, and (\ref{eq_cutoff}). Let $\tau\ll 1$, let $\{x_i:i\in I\}$ be a $(\tau r)$-net of $(X,\rho)$, and let $\{\phi_i:i\in I\}$ be the family of functions given in Lemma \ref{lem_cutpart}\ref{item_partition}. The constants implicit in the notation $\lesssim$, $\gtrsim$, and $\asymp$ below are independent of $\tau$ and $r$.

For any $f,g\in \mathcal{F}^{(p)}$, define
$$\Upsilon(f;g):=\int_W f \mathrm{d}\Gamma^{(p)}(u;g),$$
then $\Upsilon:\mathcal{F}^{(p)}\times \mathcal{F}^{(p)}\to \mathbb{R}$ is bilinear. For any $i,j\in I$, let
$$a_{ij}:=\Upsilon(\chi\phi_i;\phi_j)-\Upsilon(\chi\phi_j;\phi_i),$$
then the strongly local property yields that $a_{ij}\ne0$ only if $\rho(x_i,x_j)<4\tau r$.

For $j\in I$, consider the total variation $\lvert \Gamma^{(p)}(u;\phi_j)\rvert$ of the signed measure $\Gamma^{(p)}(u;\phi_j)$ on $W$, for any Borel subset $A\subseteq W$, by \cite[Proposition 4.8]{KS24aV5}, (\ref{eq_partition}), and (\ref{eq_mu_lower}),
\begin{align*}
&\lvert \Gamma^{(p)}(u;\phi_j)\rvert(A)\le \mu(A)^{\frac{p-1}{p}}\Gamma^{(p)}(\phi_j)(A)^{\frac{1}{p}}\lesssim \mu(A)^{\frac{p-1}{p}}\left(\frac{1}{(\tau r)^p}\mathcal{H}_\rho^1(A)\right)^{\frac{1}{p}}\lesssim \frac{1}{\tau r}\mu(A),
\end{align*}
hence
\begin{equation}\label{eq_muphi}
\lvert \Gamma^{(p)}(u;\phi_j)\rvert|_{W}\lesssim \frac{1}{\tau r}\mu|_W.
\end{equation}
Similarly,
\begin{equation}\label{eq_muchi}
\lvert \Gamma^{(p)}(u;\chi)\rvert|_{W}\lesssim \frac{1}{r}\mu|_W.
\end{equation}

For any $i,j\in I$, by the strongly local property and (\ref{eq_muphi}),
\begin{align}
&\lvert \Upsilon(\chi \phi_i;\phi_j)\rvert\nonumber\\
&\le \int_{W\cap B_\rho(x_i,2\tau r)\cap B_\rho(x_j,2\tau r)} \chi\phi_i \mathrm{d} \lvert \Gamma^{(p)}(u;\phi_j)\rvert\nonumber\\
&\lesssim \frac{1}{\tau r}\int_{W\cap B_\rho(x_i,2\tau r)\cap B_\rho(x_j,2\tau r)}\chi\phi_i \mathrm{d}\mu.\label{eq_Upij}
\end{align}

Let $(V,E)$ be a finite directed graph given by
$$V:=\left\{i\in I:B_\rho(x_i,2\tau r)\cap B_\rho(x,2r)\ne\emptyset\right\},\qquad E:=\left\{(i,j):i,j\in V,a_{ij}>0\right\}.$$
Since $(X,d)$ is doubling, (\ref{eq_rhodc}) implies that $(X,\rho)$ is also doubling. Hence there exists a positive integer $N(\tau)$, depending only on $\tau$ and the doubling constants of $(X,\rho)$, such that $\# V\le N(\tau)$. We have
\begin{equation}\label{eq_phisum}
\sum_{j\in V}\phi_j=\sum_{j\in I}\phi_j=1\text{ in }B_\rho(x,2r)\supseteq \mathrm{supp}(\chi).
\end{equation}

Let $f:E\to[0,+\infty)$ be a flow given by $f((i,j))=a_{ij}$ for $(i,j)\in E$. Then for any $j\in V$,
\begin{equation}\label{eq_parf1}
(\partial f)(j):=\sum_{i\in V}a_{ij}=\sum_{i\in V}\Upsilon(\chi\phi_i;\phi_j)-\sum_{i\in V}\Upsilon(\chi\phi_j;\phi_i)=\Upsilon(\chi;\phi_j)=\int_{W}\chi \mathrm{d}\Gamma^{(p)}(u;\phi_j).
\end{equation}
By the Leibniz rule, since $u$ is harmonic in $\Omega$ and $\chi\phi_j\in \mathcal{F}^{(p)}(\Omega)$,
\begin{equation}\label{eq_parf2}
(\partial f)(j)=\int_W \mathrm{d}\Gamma^{(p)}(u;\chi \phi_j)-\int_W\phi_j \mathrm{d}\Gamma^{(p)}(u;\chi)=-\int_{\Omega\backslash W} \mathrm{d}\Gamma^{(p)}(u;\chi \phi_j)-\int_W\phi_j \mathrm{d}\Gamma^{(p)}(u;\chi),
\end{equation}
where
$$\lvert \int_{\Omega\backslash W} \mathrm{d}\Gamma^{(p)}(u;\chi \phi_j)\rvert=\lvert \int_{B_\rho(x,2r)\backslash W} \mathrm{d}\Gamma^{(p)}(u;\chi \phi_j)\rvert\le \mu \left(B_\rho(x,2r)\backslash W\right)^{\frac{p-1}{p}} \mathcal{E}^{(p)}(\chi\phi_j)^{\frac{1}{p}}.$$
Since $0\le\chi,\phi_j\le1$ on $X$, by \ref{eq_Alg}, (\ref{eq_cutoff}), and (\ref{eq_partition}),
$$\mathcal{E}^{(p)}(\chi\phi_j)^{\frac{1}{p}}\le \mathcal{E}^{(p)}(\phi_j)^{\frac{1}{p}}+ \mathcal{E}^{(p)}(\chi)^{\frac{1}{p}}\lesssim \frac{1}{r^{\frac{p-1}{p}}}+\frac{1}{(\tau r)^{\frac{p-1}{p}}}\asymp\frac{1}{(\tau r)^{\frac{p-1}{p}}},$$
which gives
\begin{equation}\label{eq_parf3}
\sum_{j\in V}\lvert \int_{\Omega\backslash W} \mathrm{d}\Gamma^{(p)}(u;\chi \phi_j)\rvert\lesssim \frac{N(\tau)}{\tau^{\frac{p-1}{p}}}\left(\frac{\mu(B_\rho(x,2r)\backslash W)}{r}\right)^{\frac{p-1}{p}}.
\end{equation}
By (\ref{eq_phisum}), (\ref{eq_muchi}), and (\ref{eq_mu_upper}),
\begin{align*}
&\sum_{j\in V}\lvert \int_W\phi_j \mathrm{d}\Gamma^{(p)}(u;\chi)\rvert\le \sum_{j\in V}\int_W\phi_j \mathrm{d}\lvert \Gamma^{(p)}(u;\chi)\rvert\\
&=\int_{B_\rho(x,2r)\cap W}\mathrm{d}\lvert \Gamma^{(p)}(u;\chi)\rvert\lesssim \frac{1}{r}\mu(B_\rho(x,2r)\cap W)\lesssim 1.
\end{align*}
Combining this with (\ref{eq_parf2}) and (\ref{eq_parf3}),
\begin{equation}\label{eq_parf}
\sum_{j\in V} \lvert (\partial f)(j)\rvert\lesssim1+\frac{N(\tau)}{\tau^{\frac{p-1}{p}}}\left(\frac{\mu(B_\rho(x,2r)\backslash W)}{r}\right)^{\frac{p-1}{p}}.
\end{equation}

Let $l:E\to(0,+\infty)$ be a length function given by $l((i,j))=\frac{\rho(x_i,x_j)}{\tau r}\in(0,4)$ for $(i,j)\in E$. Then by (\ref{eq_Upij}), the bounded overlap property, (\ref{eq_phisum}), and (\ref{eq_mu_upper}), the transportation cost
\begin{align}
&M_l(f):=\sum_{(i,j)\in E}f((i,j))l((i,j))\lesssim \sum_{i,j\in V}\lvert \Upsilon(\chi\phi_i;\phi_j)\rvert\nonumber\\
&\lesssim \frac{1}{\tau r}\sum_{i,j\in V}\int_{W\cap B_\rho(x_i,2\tau r)\cap B_\rho(x_j,2\tau r)}\chi\phi_i \mathrm{d}\mu\lesssim\frac{1}{\tau r}\sum_{i\in V}\int_{W\cap B_\rho(x_i,2\tau r)}\chi\phi_i \mathrm{d}\mu\nonumber\\
&\lesssim \frac{1}{\tau r}\mu(B_\rho(x,2r))\lesssim\frac{1}{\tau}.\label{eq_Mlf}
\end{align}

Let $v=u-u(x)$ and let $H:V\to \mathbb{R}$ be given by $H(i)=v(x_i)$. Then by (\ref{eq_ulip}), for any $i,j\in V$,
\begin{equation}\label{eq_Hij}
\lvert H(i)-H(j)\rvert=\lvert u(x_i)-u(x_j)\rvert\lesssim \rho(x_i,x_j)\lesssim r.
\end{equation}
By the strongly local property and (\ref{eq_parf1}),
\begin{equation}\label{eq_parfH1}
(\partial f,H)_V:=\sum_{j\in V}(\partial f)(j)H(j)=\int_W\chi \mathrm{d}\Gamma^{(p)}\left(v;\sum_{j\in V}v(x_j)\phi_j\right)=\int_W\chi \mathrm{d}\Gamma^{(p)}(v;w),
\end{equation}
where
$$w:=\sum_{j\in V}v(x_j)\phi_j.$$
By (\ref{eq_muonW}), there exists a constant $C_1>0$ such that
\begin{equation}\label{eq_parfH2}
\int_W\chi \mathrm{d}\Gamma^{(p)}(v;v)\ge\mu(B_\rho(x,r)\cap W)\ge \frac{1}{C_1} r.
\end{equation}
By the Leibniz rule, since $u$ is harmonic in $\Omega$ and $\chi(v-w)\in \mathcal{F}^{(p)}(\Omega)$,
\begin{align}
&\int_W\chi \mathrm{d}\Gamma^{(p)}(v;v-w)=\int_W \mathrm{d}\Gamma^{(p)}(v;\chi(v-w))-\int_W(v-w) \mathrm{d}\Gamma^{(p)}(v;\chi)\nonumber\\
&=-\int_{\Omega\backslash W} \mathrm{d}\Gamma^{(p)}\left(v;\chi(v-w)\right)-\int_W(v-w) \mathrm{d}\Gamma^{(p)}(v;\chi).\label{eq_parfH3}
\end{align}
By (\ref{eq_phisum}) and (\ref{eq_ulip}), in $B_\rho(x,2r)$,
\begin{equation}\label{eq_parfH4}
\lvert v-w\rvert= \lvert\sum_{j\in V}(v-v(x_j))\phi_j\rvert\le \sum_{j\in V}\lvert u-u(x_j)\rvert\phi_j\lesssim \tau r,
\end{equation}
hence by (\ref{eq_muchi}) and (\ref{eq_mu_upper}),
$$\lvert \int_W(v-w) \mathrm{d}\Gamma^{(p)}(v;\chi)\rvert\le \int_{B_\rho(x,2r)\cap W} \lvert v-w\rvert \mathrm{d}\lvert\Gamma^{(p)}(v;\chi)\rvert\lesssim (\tau r) \left(\frac{1}{r}\mu(B_\rho(x,2r)\cap W)\right)\lesssim \tau r.$$
That is, there exists a constant $C_2>0$ such that
\begin{equation}\label{eq_parfH5}
\lvert \int_W(v-w) \mathrm{d}\Gamma^{(p)}(v;\chi)\rvert\le C_2\tau r.
\end{equation}
Moreover,
\begin{align}
&\lvert\int_{\Omega\backslash W} \mathrm{d}\Gamma^{(p)}\left(v;\chi(v-w)\right) \rvert=\lvert\int_{B_\rho(x,2r)\backslash W} \mathrm{d}\Gamma^{(p)}\left(v;\chi(v-w)\right) \rvert \nonumber\\
&\le \mu(B_\rho(x,2r)\backslash W)^{\frac{p-1}{p}}\Gamma^{(p)}\left(\chi(v-w)\right)(B_\rho(x,2r))^{\frac{1}{p}}.\label{eq_parfH6}
\end{align}
By \ref{eq_Alg}, (\ref{eq_parfH4}), and (\ref{eq_cutoff}),
\begin{align}
&\Gamma^{(p)}\left(\chi(v-w)\right)(B_\rho(x,2r))^{\frac{1}{p}}\nonumber\\
&\le \lVert \chi\rVert_{L^\infty(B_\rho(x,2r))}\Gamma^{(p)}\left(v-w\right)(B_\rho(x,2r))^{\frac{1}{p}}+\lVert v-w\rVert_{L^\infty(B_\rho(x,2r))}\Gamma^{(p)}\left(\chi\right)(B_\rho(x,2r))^{\frac{1}{p}}\nonumber\\
&\lesssim \Gamma^{(p)}(v-w)(B_\rho(x,2r))^{\frac{1}{p}}+\tau r^{\frac{1}{p}}.\label{eq_parfH7}
\end{align}
By (\ref{eq_parfH4}) and the bounded overlap property,
\begin{align*}
&\Gamma^{(p)}(v-w)(B_\rho(x,2r))=\Gamma^{(p)}\left(\sum_{j\in V}(v-v(x_j))\phi_j\right)(B_\rho(x,2r))\\
&\le\sum_{i\in V}\Gamma^{(p)}\left(\sum_{j:\rho(x_i,x_j)<4\tau r}(v-v(x_j))\phi_j\right)(B_\rho(x_i,2\tau r))\\
&\lesssim\sum_{i\in V}\sum_{j:\rho(x_i,x_j)<4\tau r}\Gamma^{(p)}\left((v-v(x_j))\phi_j\right)(B_\rho(x_i,2\tau r)),
\end{align*}
where for any such $i,j$, by \ref{eq_Alg}, (\ref{eq_ulip}), (\ref{eq_partition}), and (\ref{eq_mu_upper}),
\begin{align*}
&\Gamma^{(p)}\left((v-v(x_j))\phi_j\right)(B_\rho(x_i,2\tau r))^{\frac{1}{p}}\\
&\le \lVert v-v(x_j)\rVert_{L^\infty(B_\rho(x_i,2\tau r))}\Gamma^{(p)}\left(\phi_j\right)(B_\rho(x_i,2\tau r))^{\frac{1}{p}}\\
&\hspace{10pt}+\lVert \phi_j\rVert_{L^\infty(B_\rho(x_i,2\tau r))}\Gamma^{(p)}\left(v-v(x_j)\right)(B_\rho(x_i,2\tau r))^{\frac{1}{p}}\\
&\lesssim (\tau r) \frac{1}{(\tau r)^{\frac{p-1}{p}}}+(\tau r)^{\frac{1}{p}}=2(\tau r)^{\frac{1}{p}},
\end{align*}
hence
\begin{align*}
\Gamma^{(p)}(v-w)(B_\rho(x,2r))\lesssim \sum_{i\in V}\sum_{j:\rho(x_i,x_j)<4\tau r} \tau r\lesssim \sum_{i\in V}\tau r\le N(\tau)\tau r.
\end{align*}
Combining this with (\ref{eq_parfH6}) and (\ref{eq_parfH7}), we obtain a constant $C_3>0$ such that
\begin{align*}
\lvert\int_{\Omega\backslash W} \mathrm{d}\Gamma^{(p)}\left(v;\chi(v-w)\right) \rvert\le C_3 \left(\frac{\mu(B_\rho(x,2r)\backslash W)}{r}\right)^{\frac{p-1}{p}}(N(\tau)\tau)^{\frac{1}{p}}r.
\end{align*}
Together with (\ref{eq_parfH1}), (\ref{eq_parfH2}), (\ref{eq_parfH3}), and (\ref{eq_parfH5}), this yields
\begin{align}
&(\partial f,H)_V=\int_W\chi \mathrm{d}\Gamma^{(p)}(v;v)-\int_W\chi \mathrm{d}\Gamma^{(p)}(v;v-w)\nonumber\\
&\ge \left(\frac{1}{C_1}-C_2\tau -C_3 \left(\frac{\mu(B_\rho(x,2r)\backslash W)}{r}\right)^{\frac{p-1}{p}}(N(\tau)\tau)^{\frac{1}{p}}\right)r.\label{eq_parfH}
\end{align}

Fix $\tau\in(0,\frac{1}{2C_1C_2})$. Since $r=r_n\downarrow0$, (\ref{eq_muoutW}) yields an integer $N$ such that, for any $n>N$,
$$\frac{N(\tau)}{\tau^{\frac{p-1}{p}}}\left(\frac{\mu(B_\rho(x,2r)\backslash W)}{r}\right)^{\frac{p-1}{p}}\le 1,\qquad C_3 \left(\frac{\mu(B_\rho(x,2r)\backslash W)}{r}\right)^{\frac{p-1}{p}}(N(\tau)\tau)^{\frac{1}{p}}\le \frac{1}{4C_1}.$$
Then (\ref{eq_parf}) and (\ref{eq_parfH}) give
$$\sum_{j\in V} \lvert (\partial f)(j)\rvert\lesssim1,\qquad (\partial f,H)_V\ge \frac{1}{4C_1}r.$$
Combining these two inequalities with (\ref{eq_Mlf}) and (\ref{eq_Hij}) and applying Proposition \ref{prop_flow}, we obtain $i_0,\ldots,i_m\in V$ with $(i_{k-1},i_k)\in E$ for any $k=1,\ldots,m$ such that
$$\sum_{k=1}^m l(i_{k-1},i_k)\lesssim \frac{1}{\tau},\qquad H(i_m)-H(i_0)\gtrsim r,$$
that is,
$$\sum_{k=1}^m \rho(x_{i_{k-1}},x_{i_k})\lesssim r,\qquad u(x_{i_m})-u(x_{i_0})\gtrsim r.$$
It follows from (\ref{eq_ulip}) and the second inequality above that
$$\rho(x_{i_0},x_{i_m})\gtrsim r.$$
For any $k=1,\ldots,m$, since $(i_{k-1},i_k)\in E$, $\rho(x_{i_{k-1}},x_{i_k})\in(0,4\tau r)$. Then by (\ref{eq_rhodc}) and the triangle inequality for the metric $d$,
\begin{align*}
&r\lesssim\rho(x_{i_0},x_{i_m})\asymp d(x_{i_0},x_{i_m})^c \le \left(\sum_{k=1}^m d(x_{i_{k-1}},x_{i_k})\right)^c\asymp \left(\sum_{k=1}^m \rho(x_{i_{k-1}},x_{i_k})^{\frac{1}{c}}\right)^c\\
&\le \left(\max_{k=1,\ldots,m}\rho(x_{i_{k-1}},x_{i_k})^{\frac{1}{c}-1}\sum_{k=1}^m \rho(x_{i_{k-1}},x_{i_k})\right)^c\lesssim \left((\tau r)^{\frac{1}{c}-1}r\right)^c=\tau^{1-c}r,
\end{align*}
hence $1\lesssim \tau^{1-c}$ for any $\tau\in(0,\frac{1}{2C_1C_2})$. Since $c\in(0,1)$, letting $\tau\downarrow0$, we obtain a contradiction.
\end{proof}

\section*{AI statement}

The author’s study of the regularity of the $p$-walk dimension began in February 2025, motivated by a preprint version of \cite{Fen26}. In Section \ref{sec_gammap3}, the idea of ``tangent current'' was suggested by OpenAI's ChatGPT (GPT-5.6 Sol) in mid-August 2026, with reference to \cite{PS13}. The author subsequently simplified the idea to obtain the elementary formulation presented in Proposition \ref{prop_flow}. During the final preparation of this paper, ChatGPT (GPT-6 Astra) was also used to assist with proofreading, polishing the English, and generating TikZ code for figures. The author takes full responsibility for the content of this paper.

\bibliographystyle{plain}

\end{document}